\documentclass[reqno,10pt,centertags,draft]{amsart}
\usepackage{amsmath,amsthm,amscd,amssymb,latexsym,upref}

\newcommand{\bbN}{{\mathbb{N}}}
\newcommand{\bbR}{{\mathbb{R}}}

\newcommand{\bbZ}{{\mathbb{Z}}}
\newcommand{\bbC}{{\mathbb{C}}}

\newcommand{\cA}{{\mathcal A}}
\newcommand{\cB}{{\mathcal B}}

\newcommand{\cF}{{\mathcal F}}

\newcommand{\cH}{{\mathcal H}}

\newcommand{\cM}{{\mathcal M}}

\newcommand{\dott}{\,\cdot\,}
\newcommand{\no}{\notag}
\newcommand{\lb}{\label}
\newcommand{\f}{\frac}

\newcommand{\ol}{\overline}

\newcommand{\loc}{\text{\rm{loc}}}

\newcommand{\Arc}{\text{\rm{Arc}}}

\newcommand{\ess}{\text{\rm{ess}}}

\newcommand{\supp}{\text{\rm{supp}}}

\newcommand{\bi}{\bibitem}
\newcommand{\hatt}{\widehat}
\newcommand{\beq}{\begin{equation}}
\newcommand{\eeq}{\end{equation}}
\newcommand{\ba}{\begin{align}}
\newcommand{\ea}{\end{align}}

\newcommand{\tr}{\text{\rm{tr}}}

\newcommand{\abs}[1]{\lvert#1\rvert}

\renewcommand{\Re}{\text{\rm Re}}
\renewcommand{\Im}{\text{\rm Im}}
\renewcommand{\ln}{\text{\rm ln}}

\newcommand{\Om}{\Omega}
\newcommand{\si}{\sigma}
\newcommand{\om}{\omega}
\newcommand{\la}{\lambda}
\newcommand{\al}{\alpha}
\newcommand{\be}{\beta}
\newcommand{\ga}{\gamma}
\newcommand{\De}{\Delta}
\newcommand{\de}{\delta}
\newcommand{\te}{\theta}
\newcommand{\ze}{\zeta}

\newcommand{\C}{\mathbb{C}}
\newcommand{\Cz}{\C\backslash\{0\}}
\newcommand{\R}{\mathbb{R}}
\newcommand{\D}{\mathbb{D}}
\newcommand{\dD}{{\partial\hspace*{.2mm}\mathbb{D}}}
\newcommand{\Z}{\mathbb{Z}}
\newcommand{\N}{\mathbb{N}}

\newcommand{\fM}{\mathfrak{M}}

\newcommand{\UU}{\mathbb{U}}
\newcommand{\Cl}{\mathbb{C}_{\ell}}
\newcommand{\Cr}{\mathbb{C}_{r}}

\newcommand{\ltz}{{\ell^2(\Z)}}
\newcommand{\Lt}[1]{{L^2(\dD;d\mu_{#1}(\cdot,k_0))}}
\newcommand{\st}{\,|\,}

\allowdisplaybreaks
\numberwithin{equation}{section}

\newtheorem{theorem}{Theorem}[section]

\newtheorem{lemma}[theorem]{Lemma}
\newtheorem{corollary}[theorem]{Corollary}
\newtheorem{hypothesis}[theorem]{Hypothesis}
\theoremstyle{definition}
\newtheorem{definition}[theorem]{Definition}
\newtheorem{remark}[theorem]{Remark}

\begin{document}

\title[A Borg-Type Theorem for CMV Operators]{A Borg-Type Theorem
Associated with Orthogonal Polynomials on the Unit Circle}
\author[F.\ Gesztesy, and M.\ Zinchenko]{Fritz
Gesztesy and Maxim Zinchenko}
\address{Department of Mathematics,
University of Missouri, Columbia, MO 65211, USA}
\email{fritz@math.missouri.edu}
\urladdr{http://www.math.missouri.edu/personnel/faculty/gesztesyf.html}
\address{Department of Mathematics,
University of Missouri, Columbia, MO 65211, USA}
\email{maxim@math.missouri.edu}
\thanks{Based upon work supported by the US National Science
Foundation under Grant No.\ DMS-0405526.}
\date{November 15, 2004}

\begin{abstract}
We prove a general Borg-type result for reflectionless unitary CMV
operators $U$ associated with orthogonal polynomials on the unit circle.
The spectrum of $U$ is assumed to be a connected arc on the unit circle.
This extends a recent result of Simon in connection with a
periodic CMV operator with spectrum the whole unit circle.

In the course of deriving the Borg-type result we also use
exponential Herglotz representations of Caratheodory functions to prove an
infinite sequence of trace formulas connected with the CMV operator $U$.
\end{abstract}

\maketitle

\section{Introduction}\label{s1}

The aim of this paper is to prove a Borg-type uniqueness theorem for
a special class of unitary doubly infinite five-diagonal
matrices. The corresponding unitary semi-infinite five-diagonal matrices
were first introduced by Cantero, Moral, and Vel\'azquez (CMV) in
\cite{CMV03}. In \cite[Sects.\ 4.5, 10.5]{Si04}, Simon introduced the
corresponding notion of unitary doubly infinite five-diagonal matrices
and coined the term ``extended'' CMV matrices. To simplify notations we
will often just speak of CMV operators whether or not they are
half-lattice of full-lattice operators indexed by $\bbN$ or $\bbZ$,
respectively.

Before we turn to Borg-type theorems, we briefly introduce the CMV
operator $U$ studied in this paper.

We denote by $\D$ the open unit disk in $\bbC$ and let
$\alpha$ be a sequence of complex numbers in $\D$,
$\alpha=\{\al_k\}_{k \in \Z} \subset \D$. The unitary CMV
operator $U$ on $\ell^2(\bbZ)$ then can be written as a
special five-diagonal doubly infinite matrix in the
standard basis of $\ell^2(\bbZ)$ according to \cite[Sects.\
4.5, 10.5]{Si04} as
\begin{align}
U = \begin{pmatrix} \ddots &&\hspace*{-8mm}\ddots
&\hspace*{-10mm}\ddots &\hspace*{-12mm}\ddots
&\hspace*{-14mm}\ddots &&&
\raisebox{-3mm}[0mm][0mm]{\hspace*{-6mm}{\Huge $0$}}
\\
&0& -\al_{0}\rho_{-1} & -\ol{\al_{-1}}\al_{0} & -\al_{1}\rho_{0} &
\rho_{0}\rho_{1}
\\
&& \rho_{-1}\rho_{0} &\ol{\al_{-1}}\rho_{0} &
-\ol{\al_{0}}\al_{1} & \ol{\al_{0}}\rho_{1} & 0
\\
&&0& -\al_{2}\rho_{1} & -\ol{\al_{1}}\al_{2} &
-\al_{3}\rho_{2} & \rho_{2}\rho_{3}
\\
&\raisebox{-4mm}[0mm][0mm]{\hspace*{-6mm}{\Huge $0$}} &&
\rho_{1}\rho_{2} & \ol{\al_{1}}\rho_{2} & -\ol{\al_{2}}\al_{3} &
\ol{\al_{2}}\rho_{3} &0&
\\
&&&&\hspace*{-14mm}\ddots &\hspace*{-14mm}\ddots
&\hspace*{-14mm}\ddots &\hspace*{-8mm}\ddots &\ddots
\end{pmatrix}. \lb{1.1}
\end{align}
Here the sequence of positive real numbers $\{\rho_k\}_{k\in\bbZ}$
is defined by
\begin{equation}
\rho_k = \sqrt{1-\abs{\al_k}^2}, \quad k\in\bbZ, \label{1.2}
\end{equation}
and terms of the form $-\ol{\alpha_k}\alpha_{k+1}$, $k\in\Z$,
represent the diagonal entries in the infinite matrix \eqref{1.1}.

The relevance of this unitary operator $U$ in
$\ell^2(\bbZ)$, more precisely, the relevance of the
corresponding half-lattice CMV operator $U_{+,0}$ in
$\ell^2(\bbN_0)$ (cf.\ \eqref{B.17}) is derived from its
intimate relationship with the trigonometric moment problem
and hence with finite measures on the unit circle $\dD$.
(Here $\bbN_0=\bbN\cup\{0\}$.)  Let
$\{\alpha_k\}_{k\in\bbN}\subset\D$  and define the transfer
matrix
\begin{equation}
S(\zeta,k)=\begin{pmatrix} \zeta & \alpha_k \\ \overline{\alpha_k}
\zeta &1
\end{pmatrix}, \quad \zeta\in\dD,\; k\in\bbN, \lb{1.3}
\end{equation}
with spectral parameter $\zeta\in\dD$. Consider the system of difference
equations
\begin{equation}
\begin{pmatrix}\varphi_+(\zeta,k) \\ \varphi^*_+(\zeta,k)\end{pmatrix}
= S(\zeta,k)\begin{pmatrix} \varphi_+(\zeta,k-1) \\
\varphi^*_+(\zeta,k-1) \end{pmatrix}, \quad \zeta\in\dD,\;
k\in\bbN  \lb{1.4}
\end{equation}
with initial condition
\begin{equation}
\begin{pmatrix} \varphi_+(\zeta,0)\\ \varphi^*_+(\zeta,0) \end{pmatrix}
=\begin{pmatrix} 1 \\ 1 \end{pmatrix}, \quad \zeta\in\dD. \lb{1.4a}
\end{equation}
Then $\varphi_+ (\dott,k)$ are monic
polynomials of degree $k$ and
\begin{equation}
\varphi^*_+ (\zeta,k)= \zeta^k \ol{\varphi_+(1/\zeta,k)},
\quad \zeta\in\dD,\; k\in\bbN_0, \lb{1.4b}
\end{equation}
the reversed ${}^*$-polynomial of $\varphi_+(\cdot,k)$, is at
most of degree $k$. These polynomials were first introduced by Szeg\H o in the
1920's in his work on the asymptotic distribution of eigenvalues of
sections of Toeplitz forms \cite{Sz20}, \cite{Sz21}  (see
also \cite[Chs.\ 1--4]{GS84}, \cite[Ch.\ XI]{Sz78}). Szeg\H o's
point of departure was the trigonometric moment problem and hence the
theory of orthogonal polynomials on the unit circle: Given a probability
measure $d\sigma_+$  supported on an infinite set on the unit circle, find
monic polynomials of degree $k$ in $\zeta=e^{i\theta}$, $\theta\in
[0,2\pi]$, such that
\begin{equation}
\int_{0}^{2\pi} d\sigma_+(e^{i\theta}) \,
\overline{\varphi_+ (e^{i\theta},k)} \varphi_+ (e^{i\theta},k')
=\gamma_{k}^{-2} \delta_{k,k'}, \quad k,k'\in\bbN_0,
\lb{1.4c}
\end{equation}
where (cf.\ \eqref{1.2})
\begin{equation}
\gamma_k^2=\begin{cases} 1, & k=0, \\
\prod_{j=1}^k \rho_j^{-2}, & k\in\bbN.
\end{cases} \lb{1.4d}
\end{equation}
One then also infers
\begin{equation}
\int_{0}^{2\pi} d\sigma_+(e^{i\theta}) \,
\overline{\varphi^*_+ (e^{i\theta},k)} \varphi^*_+ (e^{i\theta},k')
=\gamma_{k}^{-2} \delta_{k,k'}, \quad k,k'\in\bbN_0 \lb{1.4e}
\end{equation}
and obtains that $\varphi_+(\cdot,k)$ is orthogonal to
$\{\zeta^j\}_{j=0,\dots,k-1}$ in $L^2(\dD;d\sigma_+)$ and
$\varphi^*_+(\cdot,k)$ is orthogonal to
$\{\zeta^j\}_{j=1,\dots,k}$ in $L^2(\dD;d\sigma_+)$. Additional
comments in this context will be provided in Remark \ref{rB.6}.
For a detailed account of the relationship of
$U_{+,0}$ with orthogonal polynomials on the unit circle we refer to the
monumental upcoming two-volume treatise by Simon \cite{Si04} (see also
\cite{Si04b} for a description of some of the principal results in
\cite{Si04}) and the exhaustive bibliograhy therein. For classical
results on orthogonal polynomials on the unit circle we refer, for
instance, to \cite{Ak65}, \cite{Ge46}--\cite{Ge61}, \cite{GS84},
\cite{Kr45}, \cite{Sz20}--\cite{Sz78}, \cite{To63}--\cite{Ve36}. More recent
references relevant to the spectral theoretic content of this paper are
\cite{GJ96}, \cite{GJ98},
\cite{GT94}, \cite{GZ05}, \cite{GN01}, \cite{PY04}, \cite{Si04a}. Moreover,
CMV operators are intimately related to a completely integrable version
of the defocusing nonlinear Schrodinger equation (continuous in time
but discrete in space), a special case of the Ablowitz--Ladik system.
Relevant references in this context are, for instance, \cite{AL75},
\cite{APT04}, \cite{GGH05}, \cite{GH05}, \cite{MEKL95}, \cite{NS05}, and
the literature cited therein.

We note that $S(\zeta,k)$ in \eqref{1.3} is not the
transfer matrix that leads to the half-lattice CMV operator
$U_{+,0}$ in $\ell^2(\bbN_0)$ (cf.\ \eqref{B.18}). After a
suitable change of basis introduced by Cantero, Moral, and
Vel\'azquez \cite{CMV03}, the transfer matrix $S(\zeta,k)$
turns into $T(\zeta,k)$ as defined in \eqref{A.31}.

Having introduced the notion of CMV operators, we now turn to
Borg-type uniqueness theorems. From the outset, Borg-type theorems
are inverse spectral theory assertions which typically prescribe a
connected interval (or arc) as the spectrum of a self-adjoint (or
unitary) differential or difference operator, and under a
reflectionless condition imposed on the operator (one may think of a
periodicity condition on the (potential) coefficients of the
differential or difference operator) infers the explicit form of the
coefficients of the operator in question. Typically, the form of the
coefficients determined in this context is fairly simple (and usually
given by constants or functions of exponential type).

Next, we briefly describe the history of Borg-type theorems relevant
to this paper. In 1946, Borg \cite{Bo46} proved, among a variety of
other inverse spectral theorems, the following result for one-dimensional
Schr\"odinger operators. (Throughout this paper we denote by
$\sigma(\cdot)$ and $\sigma_{\ess}(\cdot)$ the spectrum and essential
spectrum of a densely defined closed linear operator in a complex
separable Hilbert space.)

\begin{theorem}[\cite{Bo46}] \lb{t1.1} ${}$ \\
Let $q\in L^1_{\loc} (\bbR)$ be real-valued and periodic. Let
$H=-\f{d^2}{dx^2}+q$ be the associated self-adjoint
Schr\"odinger operator in $L^2(\bbR)$ and suppose that
\begin{equation}
\sigma(H)=[e_0,\infty) \, \text{ for some $e_0\in\bbR$.}
\end{equation}
Then $q$ is of the form,
\begin{equation}
q(x)=e_0 \, \text{ for a.e.\ $x\in\bbR$}. \lb{1.5}
\end{equation}
\end{theorem}

Traditionally, uniqueness results such as
Theorem\ \ref{t1.1} are called Borg-type theorems. However, this
terminology is not uniquely adopted  and hence a bit unfortunate.
Indeed, inverse spectral results on finite intervals in which the
potential coefficient(s) are recovered from two spectra, were also
pioneered by Borg in his celebrated paper \cite{Bo46}, and hence are also
coined Borg-type theorems in the literature, see, e.g.,
\cite{Ma94}, \cite{Ma99a}.

A closer examination of the proof of Theorem\ \ref{t1.1} in
\cite{CGHL00} shows that periodicity of $q$ is not the point for the
uniqueness result \eqref{1.5}. The key ingredient (besides
$\sigma(H)=[e_0,\infty)$ and $q$ real-valued) is the fact that
\begin{equation}
\text{for all $x\in \bbR$, } \, \xi(\lambda,x)=1/2 \, \text{ for a.e.\
$\lambda\in\sigma_{\ess}(h)$.} \lb{1.6}
\end{equation}
Here $\xi(\lambda,x)$, the argument of the boundary value
$g(\lambda+i0,x)$ of the diagonal Green's function of $H$ on the real
axis (where $g(z,x)=(H-zI)^{-1}(x,x))$, $z\in\bbC\backslash\sigma(h)$,
$x\in\R$), is defined by
\begin{equation}
\xi(\lambda,x)=\pi^{-1}\lim_{\varepsilon\downarrow 0}
\Im(\ln(g(\lambda+i\varepsilon,x)) \, \text{
for a.e.\ $\lambda\in\bbR$ and all $x\in\R$}. \lb{1.7}
\end{equation}

Real-valued periodic potentials are known to satisfy \eqref{1.6}, but
so do certain classes of real-valued quasi-periodic and
almost-periodic potentials $q$. In particular, the class of real-valued
algebro-geometric finite-gap KdV potentials $q$ (a subclass of the set
of real-valued quasi-periodic  potentials) is a prime example
satisfying \eqref{1.6} without necessarily being periodic.
Traditionally, potentials $q$ satisfying \eqref{1.6} are called
\textit{reflectionless} (see \cite{CG02}, \cite{CGHL00} and the
references therein).

The extension of Borg's Theorem\ \ref{t1.1} to periodic matrix-valued
Schr\"odinger operators was proved by D\'epres \cite{De95}. A new
strategy of the proof based on exponential Herglotz representations
and a trace formula (cf.\ \cite{GS96}) for such potentials, as well as
the extension to reflectionless matrix-valued potentials, was
obtained in \cite{CGHL00}.

The direct analog of Borg's Theorem\ \ref{t1.1} for periodic Jacobi
operators was proved by Flaschka \cite{Fl75} in 1975.

\begin{theorem}  [\cite{Fl75}] \lb{t1.2} ${}$ \\
Suppose $a=\{a_k\}_{k\in\bbZ}$ and $b=\{b_k\}_{k\in\bbZ}$
are  periodic real-valued sequences in
$\ell^\infty(\bbZ)$ with the same period and $a_k>0$, $k\in\bbZ$. Let
$H=aS^+ +a^-S^- +b$ be the associated self-adjoint  Jacobi operator
on $\ell^2(\bbZ)$ and suppose that
\begin{equation}
\sigma(H)=[E_-,E_+] \, \text{ for some $E_-<E_+$.}
\end{equation}
Then $a=\{a_k\}_{k\in\bbZ}$ and $b=\{b_k\}_{k\in\bbZ}$ are of the form,
\begin{equation}
a_k=(E_+-E_-)/4 , \quad  b_k=(E_-+E_+)/2, \quad k\in\bbZ.  \lb{1.8}
\end{equation}
\end{theorem}

The extension of Theorem\ \ref{t1.2} to reflectionless scalar Jacobi
operators is due to Teschl \cite[Corollary\ 6.3]{Te98} (see also
\cite[Corollary\ 8.6]{Te00}). The extension of Theorem\ \ref{t1.2} to
matrix-valued reflectionless Jacobi operators (and a corresponding
result for Dirac-type difference operators) has recently been obtained
in \cite{CGR04}.

The following very recent result of Simon is the first in connection
with orthogonal polynomials on the unit circle.

\begin{theorem} [\cite{Si04}, Sect.\ 11.14] \lb{t1.3} ${}$ \\
Suppose $\alpha=\{\alpha_k\}_{k\in\bbZ}\subset\D$ is a periodic
sequence. Let $U$ be the associated unitary
CMV operator \eqref{2.3} $($cf.\ also \eqref{A.23}$)$ on
$\ell^2(\bbZ)$ and suppose that
\begin{equation}
\sigma(U)=\dD.
\end{equation}
Then $\alpha=\{\alpha_k\}_{k\in\bbZ}$  is of the form,
\begin{equation}
\alpha_k=0, \quad k\in\bbZ. \lb{1.9}
\end{equation}
\end{theorem}

We will extend Simon's result to reflectionless Verblunsky
coefficients corresponding to a CMV operator with spectrum a
connected arc on the unit circle in our principal Section \ref{s5}.

In Section \ref{s2} we prove an infinite sequence of trace formulas
connected with CMV operators $U$ using Weyl--Titchmarsh functions (and
their  exponential Herglotz representations) associated with $U$.
Section \ref{s3} proves certain scaling results for Schur functions
associated with $U$ using a Riccati-type equation for the Verblunsky
coefficients $\alpha$. The notion of reflectionless CMV operators $U$
is introduced in Section \ref{s4} and a variety of necessary
conditions (many of them also sufficient) for $U$ to be
reflectionless are established. In our principal Section \ref{s5} we
extend Simon's Borg-type result, Theorem \ref{t1.3}, from periodic
to reflectionless Verblunsky coefficients, and then we  prove our
main new result, a Borg-type theorem for reflectionless CMV operators
whose spectrum consists of a connected subarc of the unit circle
$\dD$. Appendix \ref{A} summarizes basic facts on Caratheodory and
Schur functions relevant to this paper, and Appendix \ref{B} provides
some elements of Weyl--Titchmarsh theory for CMV operators on
half-lattices and on $\bbZ$ (a detailed treatment of this material can be
found in \cite{GZ05}).

\section{Trace Formulas} \label{s2}

In this section we discuss trace formulas associated with the CMV
operator $U$ on $\ell^2(\Z)$. We freely use the notation established
in Appendices \ref{A} and \ref{B}.

\begin{hypothesis} \lb{h2.1}
Let $\alpha$ be a sequence of complex numbers such that
\begin{equation} \label{2.1}
\alpha=\{\al_k\}_{k \in \Z} \subset \D.
\end{equation}
\end{hypothesis}

Given a sequence $\alpha$ satisfying \eqref{2.1}, we
define the sequence of positive real numbers
$\{\rho_k\}_{k\in\bbZ}$ by
\begin{equation}
\rho_k = \sqrt{1-\abs{\al_k}^2}, \quad k\in\bbZ. \label{2.2}
\end{equation}

As discussed in \eqref{A.19}--\eqref{A.23}, the unitary CMV operator
$U$ on $\ell^2(\bbZ)$ then can be written as a special five-diagonal
doubly infinite matrix in the standard basis of $\ell^2(\bbZ)$ (cf.\
\cite[Sects.\ 4.5, 10.5 ]{Si04}) as,
\begin{align}
U = \begin{pmatrix} \ddots &&\hspace*{-8mm}\ddots
&\hspace*{-10mm}\ddots &\hspace*{-12mm}\ddots
&\hspace*{-14mm}\ddots &&&
\raisebox{-3mm}[0mm][0mm]{\hspace*{-6mm}{\Huge $0$}}
\\
&0& -\al_{0}\rho_{-1} & -\ol{\al_{-1}}\al_{0} & -\al_{1}\rho_{0} &
\rho_{0}\rho_{1}
\\
&& \rho_{-1}\rho_{0} &\ol{\al_{-1}}\rho_{0} &
-\ol{\al_{0}}\al_{1} & \ol{\al_{0}}\rho_{1} & 0
\\
&&0& -\al_{2}\rho_{1} & -\ol{\al_{1}}\al_{2} &
-\al_{3}\rho_{2} & \rho_{2}\rho_{3}
\\
&\raisebox{-4mm}[0mm][0mm]{\hspace*{-6mm}{\Huge $0$}} &&
\rho_{1}\rho_{2} & \ol{\al_{1}}\rho_{2} & -\ol{\al_{2}}\al_{3} &
\ol{\al_{2}}\rho_{3} &0&
\\
&&&&\hspace*{-14mm}\ddots &\hspace*{-14mm}\ddots
&\hspace*{-14mm}\ddots &\hspace*{-8mm}\ddots &\ddots
\end{pmatrix}. \lb{2.3}
\end{align}
Here terms of the form $-\ol{\alpha_k}\alpha_{k+1}$, $k\in\Z$,
represent the diagonal entries in the infinite matrix \eqref{2.3}. The
half-lattice (i.e., semi-infinite) version of $U$ was first introduced by
Cantero, Moral, and Vel\'azquez \cite{CMV03}.

Next, we recall the half-lattice Weyl--Titchmarsh functions
$M_\pm(\cdot,k)$ associated with $U$ (cf.\
\eqref{A.68}--\eqref{A.72}) and the Caratheodory function
$M_{1,1}(\cdot,k)$ in \eqref{A.94}. By Theorem \ref{tA.3}
and the fact that $M_{1,1}(0,k)=1$ by \eqref{A.98}, one
then obtains for the exponential Herglotz representation of
$M_{1,1}(\cdot,k)$, $k\in\bbZ$,
\begin{align}
& -i\ln[iM_{1,1}(z,k)]=\oint_{\dD} d\mu_0(\zeta)\,
\Upsilon_{1,1}(\zeta,k)
\f{\zeta+z}{\zeta-z}, \quad z\in\D, \lb{2.4} \\
& \;\, 0\leq \Upsilon_{1,1}(\zeta,k) \leq \pi \, \text{ for $\mu_0$-a.e.\
$\zeta\in\dD$}. \lb{2.5}
\end{align}
For our present purpose it is more convenient to rewrite \eqref{2.4}
in the form ($k\in\bbZ$)
\begin{align}
& \ln[M_{1,1}(z,k)]=i \oint_{\dD} d\mu_0(\zeta)\,
\Xi_{1,1}(\zeta,k) \f{\zeta+z}{\zeta-z}, \quad z\in\D,
\lb{2.6} \\ & \;\, -\pi/2 \leq \Xi_{1,1}(\zeta,k) \leq
\pi/2 \, \text{ for $\mu_0$-a.e.\ $\zeta\in\dD$},  \lb{2.7}
\end{align}
where
\begin{align}
\Xi_{1,1}(\zeta,k)&=\lim_{r\uparrow 1} \Im[\ln(M_{1,1}(r\zeta)] \\
&=\Upsilon_{1,1}(\zeta,k)-(\pi/2)  \, \text{ for $\mu_0$-a.e.\
$\zeta\in\dD$}. \lb{2.8}
\end{align}
We note that $M_{1,1}(0,k)=1$ also implies
\begin{equation}
\oint_{\dD} d\mu_0(\zeta)\, \Xi_{1,1}(\zeta,k)=0, \quad k\in\bbZ.
\lb{2.9}
\end{equation}

To derive trace formulas for $U$ we now expand $M_{1,1}(z,k)$ near
$z=0$. Using \eqref{A.94} one obtains
\begin{align}
M_{1,1}(z,k) & = (\delta_k,(U+zI)(U-zI)^{-1}\delta_k)_{\ell^2(\Z)} \no \\
&= 1+ 2(\delta_k,zU^*(I-zU^*)^{-1})\delta_k)_{\ell^2(\Z)} \no \\
&= 1+\sum_{j=1}^\infty M_j(U,k)z^j, \quad z\in\D, \lb{2.10}
\end{align}
where
\begin{equation}
M_j(U,k) = 2(\delta_k,(U^*)^j\delta_k)_{\ell^2(\Z)}, \quad j\in\N,\;
k\in\bbZ
\lb{2.11}
\end{equation}
and \eqref{2.11} represents a convergent expansion in
$\cB(\ell^2(\bbZ))$. (Here $\cB(\cH)$ denotes the Banach space of
bounded linear operators mapping the Hilbert space $\cH$ into itself.)
Explicitly, one computes
\begin{equation}
M_1(U,k)=-2\alpha_k\ol{\alpha_{k+1}}, \quad k\in\bbZ. \lb{2.12}
\end{equation}
Next, we recall the well-known fact that the convergent Taylor expansion
\begin{align}
                 g(z)= & 1+ \sum_{j=1}^\infty c_j z^{j}, \quad z\in\D
\lb{2.13} \\
\intertext{implies the absolutely convergent expansion}
\ln [g(z)] =
&\sum_{j=1}^\infty d_j z^{j}, \quad |z|<\varepsilon
\lb{2.14}
\end{align}
(for $\varepsilon=\varepsilon(g)$ sufficiently small), where $d_j$ can
be recursivly computed via
\begin{equation}
d_1 = c_1, \quad d_j = c_j - \sum_{\ell=1}^{j-1}
(\ell/j) c_{j-\ell} \, d_\ell, \quad j=2,3,\dots \, . \lb{2.15}
\end{equation}
Thus, one obtains
\begin{equation}
\ln(M_{1,1}(z,k))=\sum_{j=1}^\infty L_j(U,k)z^j, \,
\text{ $|z|$ sufficiently small}, \; k\in\bbZ, \lb{2.16}
\end{equation}
where
\begin{align}
L_1(U,k) & = M_1(U,k), \lb{2.17} \\
L_j(U,k) & = M_j(U,k)-\sum_{\ell=1}^{j-1} (\ell/j)M_{j-\ell}(U,k)
L_{\ell}(U,k), \quad j=2,3,\dots, \; k\in\bbZ. \no
\end{align}

\begin{theorem} \lb{t2.2}
Let $\alpha=\{\alpha_k\}_{k\in\bbZ}\subset\D$ and $k\in\Z$. Then,
\begin{equation}
L_j(U,k)=2i\oint_{\dD} d\mu_0(\zeta)\, \Xi_{1,1}(\zeta,k)\,{\ol
\zeta}^j,
\quad j\in\bbN. \lb{2.18}
\end{equation}
In particular,
\begin{equation}
L_1(U,k)=-2\alpha_k\ol{\alpha_{k+1}}=
2i\oint_{\dD} d\mu_0(\zeta)\, \Xi_{1,1}(\zeta,k)\,{\ol \zeta}. \lb{2.19}
\end{equation}
\end{theorem}
\begin{proof}
Let $z\in\D$, $k\in\bbZ$. Since
\begin{equation}
\f{\zeta+z}{\zeta-z}=1+2\sum_{j=1}^\infty (\ol\zeta z)^j, \quad
\zeta\in\dD, \lb{2.20}
\end{equation}
\eqref{2.6} implies
\begin{equation}
\ln[M_{1,1}(z,k)]=2i\sum_{j=1}^\infty \oint_{\dD}
d\mu_0(\zeta)\, \Xi_{1,1}(\zeta,k){\ol\zeta}^j z^j, \,
\text{ $|z|$ sufficiently small}. \lb{2.21}
\end{equation}
A comparison of coefficients of $z^j$ in \eqref{2.16} and \eqref{2.21}
then proves \eqref{2.18}. \eqref{2.19} is then clear from
\eqref{2.12} and \eqref{2.17}.
\end{proof}

\section{Scaling Considerations} \label{s3}

In this section we prove some facts about the scaling behavior
of the Schur functions $\Phi_\pm$ and $\Phi_{1,1}$ and use that
to obtain spectral results for $U$. Again we freely use the notation
established in Appendices \ref{A} and \ref{B}.

Throughout this section we suppose that the sequence
$\alpha$ satisfies Hypothesis \ref{h2.1}, that is,
$\alpha=\{\alpha_k\}_{k\in\bbZ}\subset\D$.

We start by recalling the Riccati-type equation
\eqref{A.80} satisfied by $\Phi_\pm$ (cf.\ Remark
\ref{rB.17}),
\begin{equation}
\alpha_k
\Phi_\pm(z,k-1)\Phi_\pm(z,k)-\Phi_\pm(z,k-1)+z\Phi_\pm(z,k)
=\ol{\alpha_k}z, \quad z\in\bbC\backslash\dD, \; k\in\bbZ.
\lb{3.1}
\end{equation}

In the following it is convenient to indicate explicitly
the $\alpha$-dependence of $\Phi_\pm$ and $\Phi_{1,1}$ and
we will thus temporarily write $\Phi_\pm(z,k;\alpha)$ and
$\Phi_{1,1}(z,k;\alpha)$, etc.

\begin{lemma} \lb{l3.1}
Let $z\in\bbC\backslash\dD$ and $k\in\bbZ$.  Suppose
$\alpha=\{\alpha_k\}_{k\in\bbZ}\subset\D$ and assume
$\{\gamma_0,\gamma_1\}\subset\dD$. Define
$\beta=\{\gamma_0\gamma_1^k\alpha_k\}_{k\in\bbZ}$. Then,
\begin{align}
\Phi_\pm(z,k;\alpha)&=\gamma_0\gamma_1^k
\Phi_\pm(\gamma_{1}z,k;\beta), \lb{3.2} \\
\Phi_{1,1}(z,k;\alpha)&= \Phi_{1,1}(\gamma_{1}z,k;\beta).
\lb{3.3}
\end{align}
\end{lemma}
\begin{proof}
We recall that
\begin{equation}
\Phi_+(\cdot,k)\colon\D\to\D, \quad
1/\Phi_-(\cdot,k)\colon\D\to\D, \quad k\in\Z,  \lb{3.4}
\end{equation}
are analytic, with unique Taylor coefficients at $z=0$, and
hence $\Phi_\pm$ are the unique solutions of the
Riccati-type equation \eqref{3.1} satisfying \eqref{3.4}.
Since the right-hand side of \eqref{3.2} also shares the
mapping properties \eqref{3.4}, it suffices to show that
the right-hand side of \eqref{3.2} satisfies the
Riccati-type equation \eqref{3.1}. Multiplying
\begin{equation}
\beta_k\Phi_\pm(z,k-1;\beta)\Phi_\pm(z,k;\beta)-\Phi_\pm(z,k-1;\beta)
+z\Phi_\pm(z,k;\beta) -z\ol{\beta_k}=0
\end{equation}
by $\gamma_0\gamma_1^{k-1}$, one infers
\begin{align}
&\beta_k\gamma_0^{-1}\gamma_1^{-k}
\big[\gamma_0\gamma_1^{k-1}\Phi_\pm(z,k-1;\beta)\big]
\big[\gamma_0\gamma_1^{k}\Phi_\pm(z,k;\beta)\big]
-\big[\gamma_0\gamma_1^{k-1}\Phi_\pm(z,k-1;\beta)\big] \no
\\ & \quad
+z\gamma_1^{-1}\big[\gamma_0\gamma_1^k\Phi_\pm(z,k;\beta)\big]
-z\gamma_1^{-1}\big[\,\ol{\beta_k}\ol{\gamma_0^{-1}\gamma_1^{-k}}\big]=0.
\end{align}
This proves \eqref{3.2}. Since $\Phi_{1,1}=\Phi_+/\Phi_-$
by \eqref{A.101}, \eqref{3.2} implies \eqref{3.3}.
\end{proof}

Next, we also indicate the explicit $\alpha$-dependence of
$U_{\pm,k_0}$ and $U$ by $U_{\pm,k_0;\alpha}$ and
$U_\alpha$, respectively. Similarly, we write
$M_\pm(z,k;\alpha)$, $M_{\ell,\ell'}(z,k;\alpha)$,
$\ell,\ell'=0,1$, and $\cM(z,k;\alpha)$.

\begin{corollary} \lb{c3.2}
Let $k_0\in\bbZ$. Suppose
$\alpha=\{\alpha_k\}_{k\in\bbZ}\subset\D$ and assume
$\{\gamma_0,\gamma_1\}\subset\dD$. Define
$\beta=\{\gamma_0\gamma_1^k\alpha_k\}_{k\in\bbZ}$. Then,
\begin{align}
\sigma_{\rm ac}(U_{\pm,k_0;\alpha})&= \gamma_1^{-1}
\sigma_{\rm ac}(U_{\pm,k_0;\beta}), \lb{3.7}
\\
\sigma(U_{\alpha})&= \gamma_1^{-1} \sigma(U_{\beta}).
\lb{3.8}
\end{align}
Moreover, the operators $U_\al$ and $\ga_1^{-1}U_\be$ are
unitarily equivalent.
\end{corollary}
\begin{proof}
Since by \eqref{A.79}, $\pm\Re(M_\pm)>0$ is equivalent to $|\Phi_\pm^{\pm
1}|<1$ and $\pm\Re(M_\pm)>0$ is equivalent to $\pm\Re(m_\pm)>0$ by
\eqref{A.69} and \eqref{A.71}, and $\{\gamma_0,\gamma_1\}\subset\dD$,
\eqref{3.7} follows from \eqref{3.2}, \eqref{A.10}, and \eqref{A.10b}.

By \eqref{3.3}, \eqref{A.6}, \eqref{A.93}, and \eqref{A.100},
\begin{align}
d\Omega_{0,0}(\ze,k;\alpha) &=
d\Omega_{0,0}(\gamma_1\ze,k;\beta), \quad \zeta\in\dD, \; k\in\Z,
\label{3.9} \\
d\Omega_{1,1}(\ze,k;\alpha) &=
d\Omega_{1,1}(\gamma_1\ze,k;\beta), \quad \zeta\in\dD, \; k\in\Z.
\lb{3.10}
\end{align}
Applying Theorem \ref{tB.22} then proves \eqref{3.8}.

Finally, we prove the unitary equivalence of $U_\al$ and
$\ga_1^{-1}U_\be$. Fix a reference point $k\in\Z$. By
\eqref{3.3} and \eqref{A.100} one then infers
\begin{equation}
M_{1,1}(z,k;\alpha)=M_{1,1}(z,k;\beta), \quad z\in\C\backslash\dD
\lb{3.11}
\end{equation}
and hence also
\begin{equation}
M_{0,0}(z,k;\alpha)=M_{0,0}(z,k;\beta), \quad z\in\C\backslash\dD,
\lb{3.12}
\end{equation}
using \eqref{A.93}. Next, using \eqref{3.2} and \eqref{A.79} one computes
\begin{equation}
M_\pm(z,k;\al)=\f{(1+\ga_0\ga_1^k)M_\pm(\ga_1 z,k;\beta)+1
-\ga_0\ga_1^k}{(1-\ga_0\ga_1^k)M_\pm(\ga_1 z,k;\beta)
+1+\ga_0\ga_1^k}, \quad z\in\C\backslash\dD.  \lb{3.13}
\end{equation}
Insertion of \eqref{3.13} into \eqref{A.97a}--\eqref{A.97c} then yields
\begin{align}
\begin{split}
M_{0,1}(z,k;\alpha)&=\begin{cases} \ga_0\ga_1^k
M_{0,1}(\ga_1 z,k;\beta), & \text{$k$ odd,} \\
\ga_0^{-1}\ga_1^{-k} M_{0,1}(\ga_1 z,k;\beta), & \text{$k$
even,}
\end{cases} \\ M_{1,0}(z,k;\alpha)&=\begin{cases}
\ga_0^{-1}\ga_1^{-k} M_{1,0}(\ga_1 z,k;\beta), & \text{$k$
odd,} \\ \ga_0\ga_1^k M_{1,0}(\ga_1 z,k;\beta), & \text{$k$
even.}
\end{cases}
\end{split}  \lb{3.14}
\end{align}
Thus,
\begin{align}
\cM(z,k;\alpha)=
\begin{cases}
\cA_k \cM(\ga_1 z,k;\beta) \cA_k^{-1}, & k \text{ odd},
\\
\cA_k^{-1} \cM(\ga_1 z,k;\beta) \cA_k, & k \text{ even},
\end{cases} & \lb{3.15}
\\ \no
& \hspace*{-1.64cm} z\in\C\backslash\dD,
\end{align}
where
\begin{equation}
\cA_k = \begin{pmatrix} (\ga_0\ga_1^k)^{-1/2} & 0 \\ 0 &
(\ga_0\ga_1^k)^{1/2} \end{pmatrix}, \quad k\in\Z.
\end{equation}
Since $\ga_0,\ga_1\in\dD$, $\cM(z,k;\alpha)$ and $\cM(\ga_1
z,k;\beta)$ are unitarily equivalent, this implies the unitary
equivalence of $U_\alpha$ and $\ga_1^{-1} U_\beta$ by
\eqref{A.88} and Theorems \ref{tA.5} and \ref{tB.22}.
\end{proof}

\section{Reflectionless Verblunsky Coefficients} \label{s4}

In this section we discuss a variety of equivalent conditions for the
Verblunsky coefficients $\alpha$ (resp., $U$) to be reflectionless.

We denote by $M_\pm(\ze,k)$, $M_{1,1}(\ze,k)$,
$\Phi_\pm(\ze,k)$, and $\Phi_{1,1}(\ze,k)$, $\ze\in\dD$,
$k\in\Z$, etc., the radial limits to the unit circle of the
corresponding functions,
\begin{align}
M_\pm(\ze,k) &= \lim_{r \uparrow 1} M_\pm(r\ze,k), &&
M_{1,1}(\ze,k) = \lim_{r \uparrow 1} M_{1,1}(r\ze,k),
\\
\Phi_\pm(\ze,k) &= \lim_{r \uparrow 1} \Phi_\pm(r\ze,k), &&
\Phi_{1,1}(\ze,k) = \lim_{r \uparrow 1} \Phi_{1,1}(r\ze,k),
\quad \ze\in\dD,\; k\in\Z.
\end{align}
These limits are known to exist $\mu_0$-almost everywhere. The
following definition of reflectionless Verblunsky
coefficients represents the analog of reflectionless coefficients
in Schr\"odinger, Dirac, and Jacobi operators (cf., e.g.
\cite{CG02}, \cite{CGHL00}, \cite{Cr89}, \cite{DS83},
\cite{GKT96}, \cite{GS96}, \cite{GJ84}, \cite{GJ86},
\cite{Jo82}, \cite{Ko84}--\cite{KK88}, \cite{SY95},
\cite{SY96a}, \cite{Te98}, \cite{Te00}).

\begin{definition} \lb{d4.1}
Let $\alpha=\{\alpha_k\}_{k\in\bbZ}\subset\D$ and denote by $U$ the
associated unitary CMV operator $U$ on $\ell^2(\Z)$. Then $\alpha$
(resp., $U$) is called {\it reflectionless}, if
\begin{equation}
\text{for all $k\in\bbZ$, } \, M_+(\zeta,k)=-\ol{M_-(\zeta,k)} \,
\text{ for $\mu_0$-a.e.\ $\zeta\in \sigma_{\ess}(U)$.} \lb{4.3}
\end{equation}
\end{definition}

The following result provides a variety of equivalent
criteria for $\alpha$ (resp., $U$) to be reflectionless.

\begin{theorem} \lb{t4.3}
Let $\alpha=\{\alpha_k\}_{k\in\bbZ}\subset\D$ and denote by
$U$ the associated unitary CMV operator $U$ on
$\ell^2(\Z)$. Then the following assertions $(i)$--$(vi)$
are equivalent:
\\
$(i)$ $\alpha=\{\alpha_k\}_{k\in\bbZ}$ is reflectionless.
\\
$(ii)$ Let $\gamma\in\dD$. Then
$\beta=\{\gamma\alpha_k\}_{k\in\bbZ}$ is reflectionless.
\\
$(iii)$ For all $k\in\bbZ$,
$M_+(\zeta,k)=-\ol{M_-(\zeta,k)}$ for $\mu_0$-a.e.\ $\zeta\in
\sigma_{\ess}(U)$.
\\
$(iv)$ For some $k_0\in\bbZ$,
$M_+(\zeta,k_0)=-\ol{M_-(\zeta,k_0)}$ for $\mu_0$-a.e.\ $\zeta\in
\sigma_{\ess}(U)$.
\\
$(v)$ For all $k\in\bbZ$,
$\Phi_+(\zeta,k)=1/\ol{\Phi_-(\zeta,k)}$ for $\mu_0$-a.e.\
$\zeta\in \sigma_{\ess}(U)$.
\\
$(vi)$ For some $k_0\in\bbZ$,
$\Phi_+(\zeta,k_0)=1/\ol{\Phi_-(\zeta,k_0)}$ for $\mu_0$-a.e.\
$\zeta\in\sigma_{\ess}(U)$.  \\[1mm]
Moreover, conditions $(i)$--$(vi)$ imply the following equivalent
assertions $(vii)$--$(ix)$:
\\
$(vii)$ For all $k\in\bbZ$, $\Xi_{1,1}(\zeta,k)=0$ for
$\mu_0$-a.e.\ $\zeta\in \sigma_{\ess}(U)$.
\\
$(viii)$ For all $k\in\bbZ$, $M_{1,1}(\zeta,k) > 0$ for
$\mu_0$-a.e.\ $\zeta\in \sigma_{\ess}(U)$.
\\
$(ix)$ For all $k\in\bbZ$, $\Phi_{1,1}(\zeta,k)\in (-1,1)$ for
$\mu_0$-a.e.\ $\zeta\in \sigma_{\ess}(U)$.
\end{theorem}
\begin{proof}
We will prove the following diagram:
\begin{align*}
&(ii)
\\
& \hspace*{1pt} \Updownarrow
\\ \no
(i) \Leftrightarrow (iii) \Leftrightarrow &(v)
\Leftrightarrow (vi) \Leftrightarrow (iv)
\\ \no
& \hspace*{1pt} \Downarrow
\\
& (ix) \Leftrightarrow (viii) \Leftrightarrow (vii)
\end{align*}

\noindent $(i)$ is equivalent to $(iii)$ by Definition \ref{d4.1}.

\noindent $(iii)$ is equivalent to $(v)$ and $(vi)$ is equivalent to
$(iv)$ by \eqref{A.78} and \eqref{A.79}.

\noindent $(v) \Leftrightarrow (ii)$: By Lemma \ref{l3.1},
\begin{align}
\Phi_+(z,k;\al)\ol{\Phi_-(z,k;\al)} =
\Phi_+(z,k;\be)\ol{\Phi_-(z,k;\be)}, \quad z\in\C, \; k\in\Z,
\end{align}
hence the fact that $(i)$ is equivalent to $(v)$ implies that
$(v)$ is equivalent to $(ii)$.

\noindent That $(v)$ implies $(vi)$ is clear.

\noindent $(vi) \Rightarrow (v)$: By \eqref{A.80},
\begin{align}
\Phi_\pm(z,k+1) &= \frac{z\ol{\al_{k+1}} +
\Phi_\pm(z,k)}{\al_{k+1} \Phi_\pm(z,k) + z}, \label{4.4}
\\
\Phi_\pm(z,k-1) &= \frac{z\ol{\al_k}
-z\Phi_\pm(z,k)}{\al_k\Phi_\pm(z,k) -1}, \quad
z\in\C\backslash\dD, \; k\in\Z. \label{4.5}
\end{align}
Taking into account $(vi)$ at the point $k_0\in\Z$,
\begin{align}
\Phi_+(\ze,k_0)\ol{\Phi_-(\ze,k_0)} = 1 \, \text{ for
$\mu_0$-a.e.\ $\ze \in\si_{\text{\rm ess}}(U)$,}
\end{align}
one proves $(vi)$ at the points $k_0\pm 1$ as follows:
\begin{align}
& \Phi_+(\ze,k_0+1)\ol{\Phi_-(\ze,k_0+1)} =
\frac{\ze\ol{\al_{k_0+1}} + \Phi_+(\ze,k_0)}{\al_{k_0+1}
\Phi_+(\ze,k_0) + \ze} \; \frac{\ol{\ze}\al_{k_0+1} +
\ol{\Phi_-(\ze,k_0)}}{\ol{\al_{k_0+1}}\,\ol{\Phi_-(\ze,k_0)} +
\ol{\ze}} \no \\
& \quad = \frac{1 + \abs{\al_{k_0+1}}^2 +
\ol{\ze}\al_{k_0+1}\Phi_+(\ze,k_0) +
\ze\ol{\al_{k_0+1}}\,\ol{\Phi_-(\ze,k_0)}}
{\abs{\al_{k_0+1}}^2 + 1 +
\ol{\ze}\al_{k_0+1}\Phi_+(\ze,k_0) +
\ze\ol{\al_{k_0+1}}\,\ol{\Phi_-(\ze,k_0)}} \no \\ & \quad =
1 \text{ for $\mu_0$-a.e.\ $\ze \in \si_{\text{\rm
ess}}(U)$,}
\\ & \Phi_+(\ze,k_0-1)\ol{\Phi_-(\ze,k_0-1)} = \abs{\ze}^2
\frac{\ol{\al_{k_0}} -
\Phi_+(\ze,k_0)}{\al_{k_0}\Phi_+(\ze,k_0) - 1} \;
{\frac{\al_{k_0} -
\ol{\Phi_-(\ze,k_0)}}{\ol{\al_{k_0}}\,\ol{\Phi_-(\ze,k_0)}
- 1}} \no \\
& \quad = \frac{\abs{\al_{k_0}}^2 + 1 -
\al_{k_0}\Phi_+(\ze,k_0) -
\ol{\al_{k_0}}\,\ol{\Phi_-(\ze,k_0)}} {\abs{\al_{k_0}}^2 +
1 - \al_{k_0}\Phi_+(\ze,k_0) -
\ol{\al_{k_0}}\,\ol{\Phi_-(\ze,k_0)}} = 1 \, \text{ for
$\mu_0$-a.e.\ $\ze \in \si_{\text{\rm ess}}(U)$.}
\end{align}
Iterating this procedure implies $(v)$.

\noindent $(ix)$ is equivalent to $(viii)$ by \eqref{A.99} and
\eqref{A.100}.

\noindent $(viii)$ is equivalent to $(vii)$ by \eqref{2.8}.

\noindent $(v) \Rightarrow (ix)$: By \eqref{A.101},
\begin{align}
\Phi_{1,1}(z,k) = \frac{\Phi_+(z,k)}{\Phi_-(z,k)} =
\frac{\Phi_+(z,k)\ol{\Phi_-(z,k)}}{\abs{\Phi_-(z,k)}^2}, 
\quad z\in\C\backslash\dD, \; k\in\Z, 
\end{align}
and hence $(v)$ implies $(ix)$.
\end{proof}

The next result shows (the expected fact) that periodic Verblunsky
coefficients are  reflectionless.

\begin{lemma} \label{l4.4}
Let $\al=\{\al_k\}_{k\in\Z}$ be a sequence of periodic
Verblunsky coefficients. Then $\al$ is reflectionless.
$($This applies, in particular, to $\al=0$.$)$
\end{lemma}
\begin{proof}
Let $\om\in\N$ denote the period of $\al=\{\al_k\}_{k\in\Z}$. Without loss
of generality we may assume $\om$ to be even. (If $\om$ is odd, we can
consider the even period $2\om$.) Then,
\begin{align} \label{4.14}
\fM(z,k_0) =
\begin{pmatrix}
\fM_{1,1}(z,k_0) & \fM_{1,2}(z,k_0)
\\
\fM_{2,1}(z,k_0) & \fM_{2,2}(z,k_0)
\end{pmatrix} = \prod_{k=1}^{\om} T(z,k_0+k), \quad
z\in\Cz,\; k_0\in\Z
\end{align}
represents the monodromy matrix of the CMV operator $U$ associated with
the sequence $\al$. By $\De(z)$ we denote the  corresponding Floquet
discriminant,
\begin{align}
\De(z) = \frac{1}{2}\tr(\fM(z,k_0)), \quad z\in\Cz.
\lb{4.15}
\end{align}
We note that $\De(z)$ does not depend on $k_0$. By
\eqref{A.31} and \eqref{4.14},
\begin{align}
\fM_{1,1}(\ze,k_0) &= \ol{\fM_{2,2}(\ze,k_0)}, \label{4.16}
\\
\fM_{1,2}(\ze,k_0) &= \ol{\fM_{2,1}(\ze,k_0)}, \quad
\ze\in\dD,\;  k_0\in\Z. \label{4.17}
\end{align}
Thus, $\De(\ze)=\Re(\fM_{1,1}(\ze,k_0))\in\R$ for all
$\ze\in\dD$. Moreover, since $\det(\fM(z,k_0))=1$, for all
$k_0\in\Z$, the eigenvalues of $\fM(z,k_0)$ are given by
\begin{align}
\rho_\pm(z) = \De(z) \mp \sqrt{\De(z)^2-1}, \quad z\in\Cz,
\end{align}
where the branch of the square root is chosen such that
$\abs{\rho_\pm(z)} \lessgtr 1$ for $z\in\bbC\backslash(\dD\cup\{0\})$, and
hence,
\begin{align}
\binom{u_\pm(z,k+\om,k_0)}{v_\pm(z,k+\om,k_0)}=\rho_\pm(z)
\binom{u_\pm(z,k,k_0)}{v_\pm(z,k,k_0)}, \quad z\in\Cz, \;
k,k_0\in\Z.
\end{align}
Thus, $\rho_\pm$ are the Floquet multipliers associated with $U$ and
consequently one obtains the following characterization of the
spectrum of $U$,
\begin{align} \label{4.18}
\si(U) = \{\ze\in\dD \st \abs{\rho_\pm(\ze)} = 1\} =
\{\ze\in\dD \st -1 \leq \De(\ze) \leq 1\}.
\end{align}
Next, assume $k$ to be even. Then \eqref{B.72} implies
\begin{align}
\Phi_\pm(z,k) &= \frac{u_\pm(z,k,k_0)}{v_\pm(z,k,k_0)} =
\frac{u_\pm(z,k+\om,k_0)/\rho_\pm(z)}{v_\pm(z,k+\om,k_0)/\rho_\pm(z)}
\no \\
&= \frac{\fM_{1,1}(z,k)\Phi_\pm(z,k)+\fM_{1,2}(z,k)}
{\fM_{2,1}(z,k)\Phi_\pm(z,k)+\fM_{2,2}(z,k)}, \quad z\in\Cz.
\end{align}
It follows that
\begin{align}
\Phi_\pm(z,k) = \frac{\fM_{1,1}(z,k)-\fM_{2,2}(z,k)
\pm2\sqrt{\De(z)^2-1}}{2\,\fM_{2,1}(z,k)}, \quad z\in\Cz,
\end{align}
and hence, by \eqref{4.16} and \eqref{4.18},
\begin{align}
\Phi_\pm(\ze,k) = i\frac{\Im(\fM_{1,1}(\ze,k))
\pm\sqrt{1-\Re(\fM_{1,1}(\ze,k))^2}}{\fM_{2,1}(\ze,k)},
\quad \ze\in\si(U).
\end{align}
Thus, by \eqref{4.16}, \eqref{4.17}, and $\det(\fM(z,k))=1$
for all $z\in\Cz$, $k\in\Z$,
\begin{align}
\Phi_+(\ze,k)\ol{\Phi_-(\ze,k)} =
\frac{\abs{\fM_{1,1}(\ze,k)}^2 -
1}{\abs{\fM_{2,1}(\ze,k)}^2} = 1, \quad \ze\in\si(U),
\lb{4.22}
\end{align}
and hence $\al$ is reflectionless by Theorem
\ref{t4.3}\,$(vi)$.
\end{proof}

We conclude this section with another result concerning the
reflectionless condition \eqref{4.3} on arcs of the unit circle. It
is contained in Lemma 10.11.17 in \cite{Si04}. The latter is based on
results in \cite{Ko84} (see also \cite{Gr60}, \cite{Ko87}). For
completeness we include the following elementary proof (which only
slightly differs from that in \cite{Si04} in that no
$H^p$-arguments are involved). To fix some notation we denote  by
$f_+$ and $f_-$ a Caratheodory and anti-Caratheodory function,
respectivley, and by
$\varphi_+$ and
$\varphi_-$ the corresponding Schur and anti-Schur function,
\begin{equation}
\varphi_\pm = \f{f_\pm-1}{f_\pm+1}.  \lb{4.26}
\end{equation}
Moreover, we introduce the corresponding Herglotz representations
of $f_\pm$ (cf.\ \eqref{A.3}, \eqref{A.4})
\begin{equation}
f_\pm(z)=ic_\pm \pm \oint_{\dD} d\mu_\pm (\zeta) \,
\f{\zeta+z}{\zeta-z}, \quad z\in\D, \; c_\pm\in\bbR. \lb{4.26a}
\end{equation}
We introduce the following notation for open arcs on the unit circle
$\dD$,
\begin{equation}
\Arc\big(\big(e^{i\theta_1},e^{i\theta_2}\big)\big)
=\big\{e^{i\theta}\in\dD\,|\,
\theta_1<\theta< \theta_2\big\}, \quad \theta_1 \in
[0,2\pi), \; \theta_1< \theta_2\leq \theta_1+2\pi. \lb{4.26b}
\end{equation}
An open arc $A\subseteq \dD$ then either coincides with
$\Arc\big(\big(e^{i\theta_1},e^{i\theta_2}\big)\big)$ for some
$\theta_1 \in [0,2\pi)$, $\theta_1< \theta_2\leq \theta_1+2\pi$, or
else, $A=\dD$.

\begin{lemma} \lb{l4.5}
Let $A\subseteq\dD$ be an open arc and assume that  $f_+$ $($resp.,
$f_-$$)$ is a Caratheodory $($resp., anti-Caratheodory$)$ function
satisfying the reflectionless condition \eqref{4.3} $\mu_0$-a.e. on
$A$, that is,
\begin{equation}
\lim_{r\uparrow 1}\big[f_+(r\zeta)+\ol{f_-(r\zeta)} \big]=0 \,
\text{ $\mu_0$-a.e.\ on $A$}.  \lb{4.27}
\end{equation}
Then, \\
(i) $f_+(\zeta)=-\ol{f_-(\zeta)}$ for all $\zeta\in A$. \\
(ii) For $z\in\D$, $-\ol{f_-(1/\ol z)}$ is the analytic
continuation of $f_+(z)$ through the arc $A$. \\
(iii) $d\mu_\pm$ are purely absolutely continuous on $A$ and
\begin{equation}
\f{d\mu_{\pm}}{d\mu_0}(\zeta)=\Re(f_+(\zeta)) =-\Re(f_-(\zeta)), \quad
\zeta\in A. \lb{4.30}
\end{equation}
\end{lemma}
\begin{proof}
By \eqref{4.26} and
\begin{equation}
\varphi_+(z) -\ol{1/\varphi_-(z)}
=\f{-2\big[f_+(z)+\ol{f_-(z)}\,\big]}{\big[f_+(z)+1\big]
\big[\,\ol{f_-(z)}-1\big]}, \quad  z\in\D,   \lb{4.31}
\end{equation}
equation \eqref{4.27} is equivalent to
\begin{equation}
\lim_{r\uparrow 1}\big[\varphi_+(r\zeta)
-\ol{1/\varphi_-(r\zeta)}\,\big]=0 \, \text{ for
$\mu_0$-a.e.\ $\zeta\in A$}.  \lb{4.32}
\end{equation}
Next, introducing
\begin{equation}
g_1(z)=[\varphi_+(z)-1/\varphi_-(z)]/2, \quad
g_2(z)=[\varphi_+(z)+1/\varphi_-(z)]/(2i), \quad z\in\D, \lb{4.33}
\end{equation}
then $g_j$, $j=1,2$, are Schur functions (since $z_1,z_2\in\D$
implies $(z_1\pm z_2)/2 \in\D$) and hence,
\begin{equation}
\text{$g_j+1$, $j=1,2$, are Caratheodory functions.}  \lb{4.34}
\end{equation}
Moreover, by
\eqref{4.32},
\begin{equation}
\Re(g_j(\zeta))=\lim_{r\uparrow 1} \Re(g_j(r\zeta))=0 \,
\text{ for $\mu_0$-a.e. }\zeta\in A, \; j=1,2. \lb{4.35}
\end{equation}
Since $g_j$, $j=1,2$, are Schur functions,
\begin{equation}
|g_j(z)|\leq 1, \quad z\in \D, \; j=1,2,  \lb{4.36}
\end{equation}
and hence the measures $d\mu_j$ in the Herglotz
representation of $g_j+1$, $j=1,2$, are purely absolutely
continuous by \eqref{A.10d},
\begin{equation}
d\mu_j=d\mu_{j,\rm ac}, \quad d\mu_{j,\rm s}=0, \quad
j=1,2. \lb{4.38}
\end{equation}
By \eqref{A.5a} and \eqref{A.10} one thus obtains
\begin{align}
g_j(z)+1&=ic_j +\oint_{\dD} [\Re(g_j(\zeta))+1]d\mu_0(\zeta)\,
\f{\zeta+z}{\zeta-z}, \no \\
&=ic_j +1 +\oint_{\dD} \Re(g_j(\zeta))d\mu_0(\zeta)\,
\f{\zeta+z}{\zeta-z}, \quad z\in\D, \; j=1,2, \lb{4.39}
\end{align}
that is,
\begin{equation}
g_j(z)=ic_j +\oint_{\dD} \Re(g_j(\zeta))d\mu_0(\zeta)\,
\f{\zeta+z}{\zeta-z}, \quad z\in\D, \; j=1,2. \lb{4.40}
\end{equation}
By \eqref{4.35}, the signed measure $\Re(g_j)d\mu_0$ has no
support on the arc $A$ and hence $g_j$, $j=1,2$, admit an analytic
continuation through $A$. Moreover, using \eqref{4.40} one computes
\begin{equation}
\ol{g_j(\zeta_0)}=-g_j(\zeta_0), \quad \zeta_0\in A, \;
j=1,2. \lb{4.41}
\end{equation}
Thus, the Schwarz symmetry principle yields
\begin{equation}
g_j(z)=-\ol{g_j(1/\ol z)}, \quad z\in\bbC\backslash\D, \; j=1,2.
\lb{4.42}
\end{equation}
Since
\begin{equation}
\varphi_+=g_1+ig_2, \quad 1/\varphi_-=-g_1+ig_2, \lb{4.43}
\end{equation}
also $\varphi_+$ and $1/\varphi_-$ admit analytic
continuations through the open arc $A$ and because of
\eqref{4.42} (and in agreement with \eqref{4.32}) one
obtains
\begin{equation}
\varphi_+(z)=\ol{1/\varphi_-(1/\ol z)}, \quad z\in\bbC\backslash\D.
\lb{4.44}
\end{equation}
Thus, one computes
\begin{equation}
f_+(z)=\f{1+\varphi_+(z)}{1-\varphi_+(z)}=
\f{1+\ol{1/\varphi_-(1/\ol z)}}{1-\ol{1/\varphi_-(1/\ol z))}}
=\f{\ol{\varphi_-(1/\ol z)}+1}{\ol{\varphi_-(1/\ol z)}-1}
=-\ol{f_-(1/\ol z)}, \quad z\in\bbC\backslash\D.  \lb{4.45}
\end{equation}
This proves items $(i)$ and $(ii)$. In particular,
$\Re(f_\pm(\zeta))$ exists and is finite for all $\zeta\in
A$ and hence
\begin{equation}
S_{\mu_{\pm,\rm s}} \cap A =\emptyset, \lb{4.46}
\end{equation}
where $S_{\mu_{\pm,\rm s}}$ denotes an essential support of
$d\mu_{\pm,\rm s}$. By \eqref{A.10} one thus computes
\begin{align}
\f{d\mu_{\pm}}{d\mu_0}(\zeta)=\Re(f_+(\zeta))=
-\Re(f_-(\zeta)), \quad \zeta\in A,  \lb{4.47}
\end{align}
proving item $(iii)$.
\end{proof}

It is perhaps worth noting that this proof is based on the
elementary fact that if $g$ is any Schur function, then
$g +1$ is a Caratheodory function with purely absolutely
continuous measure in its Herglotz representation (cf.\
the first line of \eqref{4.39}). (In particular, the support of the
measure in the Herglotz representation of $g+1$ equals $\dD$.) The
rest are simple Schwarz symmetry considerations.

\section{The Borg-type Theorem for CMV operators} \label{s5}

We recall our notation for closed arcs on the unit circle $\dD$,
\begin{equation}
\Arc\big(\big[e^{i\theta_1},e^{i\theta_2}\big]\big)
=\big\{e^{i\theta}\in\dD\,|\,
\theta_1\leq\theta\leq \theta_2\big\}, \quad \theta_1 \in
[0,2\pi), \; \theta_1\leq \theta_2\leq \theta_1+2\pi \lb{5.1}
\end{equation}
and similarly for open arcs (cf.\ \eqref{4.26b}) and arcs
open or closed at one endpoint (cf.\ \eqref{A.5}).

We start with a short proof of a recent result of Simon \cite{Si04} in
the case where $\alpha$ is a periodic sequence of Verblunsky
coefficients, see Theorem \ref{t1.3}. We will extend this result from the
periodic to the reflectionless case.

\begin{theorem}  \lb{t5.1}
Let $\alpha=\{\alpha_k\}_{k\in\bbZ}\subset\D$ be a
reflectionless sequence of Verblunsky coefficients. Let $U$
be the associated unitary CMV operator \eqref{2.3} $($cf.\
also \eqref{A.19}--\eqref{A.22}$)$ on $\ell^2(\bbZ)$ and
suppose that
\begin{equation}
\sigma(U)=\dD. \lb{5.2}
\end{equation}
Then $\alpha=\{\alpha_k\}_{k\in\bbZ}$ is of the form,
\begin{equation}
\alpha_k=0, \quad k\in\bbZ. \lb{5.3}
\end{equation}
\end{theorem}
\begin{proof}
Since by hypothesis $U$ is reflectionless, one infers from Definition
\ref{d4.1} that
\begin{equation}
\text{for all $k\in\Z$, } \, M_+(\zeta,k)=-\ol{M_-(\zeta,k)} \,
\text{ for $\mu_0$-a.e.\ $\zeta\in\dD$.} \lb{5.4}
\end{equation}
Denote by $d\omega_\pm(\cdot,k)$ the measures associated with
the Herglotz representation \eqref{A.3} of $M_\pm(\cdot,k)$,
$k\in\Z$. (Of course, $d\omega_+=d\mu_+$ by
\eqref{A.69}.) By Lemma \ref{l4.5}, $d\omega_\pm(\cdot,k)$ are
purely absolutely continuous for all $k\in\bbZ$,
\begin{equation}
d\omega_\pm(\cdot,k)=d\omega_{\pm,\rm ac}(\cdot,k), \quad
k\in\Z. \lb{5.5}
\end{equation}
Moreover, by \eqref{5.4}, \eqref{5.5}, and \eqref{A.10} one
concludes that
\begin{equation}
d\omega_+(\cdot,k)=d\omega_-(\cdot,k), \quad k\in\Z
\lb{5.6}
\end{equation}
and hence that
\begin{equation}
M_+(z,k)=-M_-(z,k), \quad z\in\C, \; k\in\Z. \lb{5.7}
\end{equation}
Taking $z=0$ in \eqref{5.7}, and utilizing \eqref{A.70} and
\eqref{A.72} then proves
\begin{equation}
1=- \f{\alpha_k+1}{\alpha_k-1}, \quad k\in\Z \lb{5.8}
\end{equation}
and hence \eqref{5.3} holds.
\end{proof}

Actually, still assuming the hypotheses of Theorem \ref{t5.1}, one 
can go a bit further: 
In addition to \eqref{5.3} and \eqref{5.4}, \eqref{5.2} and
Theorem \ref{t4.3}\,$(vii)$ yield that
\begin{equation}
\text{for all $k\in\Z$, } \, \Xi_{1,1}(\zeta,k)=0 \, \text{
for $\mu_0$-a.e.\ $\zeta\in\dD$}  \lb{5.9}
\end{equation}
and hence that
\begin{equation}
M_{1,1}(z,k)=1, \quad z\in\D, \; k\in\Z.  \lb{5.10}
\end{equation}
Moreover, \eqref{5.2}, \eqref{5.4}, and \eqref{A.97} yield
\begin{equation}
\text{for all $k\in\Z$, } \, M_\pm(z,k)=\pm 1, \; z\in\D,
\quad d\omega_\pm(\cdot,k)=d\mu_0.  \lb{5.11}
\end{equation}

\begin{remark}  \lb{r5.2}
The special case where $\alpha$ is periodic and $\sigma(U)=\dD$
and thus $\alpha=0$ has originally been derived by Simon
\cite[Sect.\ 11.14]{Si04} using different techniques based on Floquet
theory (cf.\ Theorem \ref{t1.3}).
\end{remark}

The principal result of this paper then reads as follows.

\begin{theorem}  \lb{t5.3}
Let $\alpha=\{\alpha_k\}_{k\in\bbZ}\subset\D$ be a
reflectionless sequence of Verblunsky coefficients. Let $U$
be the associated unitary CMV operator \eqref{2.3} $($cf.\
also \eqref{A.19}--\eqref{A.22}$)$ on $\ell^2(\bbZ)$ and
suppose that the spectrum of $U$ consists of a connected arc 
of $\dD$,
\begin{equation}
\sigma(U)=\Arc\big(\big[e^{i\theta_0},e^{i\theta_1}\big]\big) 
 \lb{5.12}
\end{equation}
with $\theta_0 \in [0,2\pi]$,
$\theta_0<\theta_1\leq\theta_0+2\pi$, and hence
$e^{i(\theta_0+\theta_1)/2}\in\Arc\big(\big(e^{i\theta_0},
e^{i\theta_1}\big)\big)$. Then
$\alpha=\{\alpha_k\}_{k\in\bbZ}$ is of the form,
\begin{equation}
\alpha_k=\alpha_0 g^k, \quad k\in\bbZ, \lb{5.13}
\end{equation}
where
\begin{equation}
g=-\exp(i(\theta_0+\theta_1)/2) \, \text{ and } \,
|\alpha_0|=\cos((\theta_1-\theta_0)/4). \lb{5.14}
\end{equation}
\end{theorem}
\begin{proof}
By Theorem \ref{t4.3}\,$(vii)$ (as a consequence of the
reflectionless property of $\alpha$) and the fact that
$M_{1,1}(\cdot,k)$, $k\in\Z$, is purely imaginary on the
spectral gap $\Arc((e^{i\te_1},e^{i\te_0+2\pi}))$ (since by
Theorem \ref{tB.22} $\supp(d\Om_{1,1}) \subseteq \si(U)$)
and strictly monotone as described in \ref{A.9a}, there
exists a $\theta_*(k)\in[\te_1,\te_0+2\pi]$ such that
$\Xi_{1,1}(\cdot,k)$, $k\in\Z$, is of the form
\begin{equation}
\Xi_{1,1}(\zeta,k)=\begin{cases} 0, &
\zeta\in\Arc\big(\big(e^{i\theta_0},e^{i\theta_1}\big)\big), \\
\pi/2, &
\zeta\in\Arc\big(\big(e^{i\theta_1},e^{i\theta_*(k)}\big)\big), \\
-\pi/2, &
\zeta\in\Arc\big(\big(e^{i\theta_*(k)},e^{i(\theta_0+2\pi)}\big)\big)
\end{cases} 
\lb{5.15}
\end{equation}
for $\mu_0\text{-a.e. } \ze\in\dD$, $k\in\Z$.
Taking into account \eqref{2.9} then yields
\begin{align}
0&=\oint_{\dD} d\mu_0(\zeta)\,\Xi_{1,1}(\zeta,k)= \f{1}{4}
\oint_{\theta_1}^{\theta_*(k)} d\theta - \f{1}{4}
\oint_{\theta_*(k)}^{\theta_0+2\pi} d\theta \no \\
&=\f{1}{4}[2\theta_*(k)-\theta_0-2\pi-\theta_1], \quad
k\in\Z \lb{5.16}
\end{align}
and hence
\begin{equation}
\theta_*(k)=\f{1}{2}(\theta_0+\theta_1)+\pi, \quad k\in\Z
\end{equation}
is in fact $k$-independent and denoted by $\theta_*$ in the
following. As a result, 
$\Xi_{1,1}(\cdot,k)=\Xi_{1,1}(\cdot)$ in \eqref{5.15} is
also $k$-independent.

By \eqref{2.19},
\begin{align}
\alpha_k\ol{\alpha_{k+1}} &= -i\oint_{\dD} d\mu_0(\zeta)\,
\Xi_{1,1}(\zeta)\,{\ol \zeta} = -i
\oint_{\te_1}^{\te_*}\frac{\pi}{2}e^{-it}\frac{dt}{2\pi} +
i\oint_{\te_*}^{\te_0+2\pi}\frac{\pi}{2}e^{-it}\frac{dt}{2\pi}
\no \\
&= -\frac{1}{4}e^{-i(\te_0+\te_1)/2}
\big(2+2\cos((\te_1-\te_0)/2)\big) \no \\ &=
-e^{-i(\te_0+\te_1)/2} \cos^2((\te_1-\te_0)/4), \quad
k\in\Z. \label{5.21}
\end{align}
Thus, $\alpha_{k_0}\ol{\alpha_{k_0+1}} = 0$ for some $k_0\in\bbZ$ is
equivalent to $\te_1=\te_0+2\pi$ and hence the assertions \eqref{5.13}
and \eqref{5.14} reduce to $\alpha=0$ as in Theorem \ref{t5.1}. (This is
of course consistent with \eqref{5.13}, \eqref{5.14}, since $|\alpha_0|=0$
in this case.) In the case $\alpha_k\ol{\alpha_{k+1}} \neq 0$ for all
$k\in\bbZ$, it follows from \eqref{5.21} that,
\begin{align*}
\al_k = \ga_0\ga_1^k \quad \text{and} \quad \cdots =
\abs{\al_1}\abs{\al_2} = \abs{\al_2}\abs{\al_3} =
\abs{\al_3}\abs{\al_4} = \abs{\al_4}\abs{\al_5} = \cdots
\end{align*}
Hence, 
\begin{align} \label{5.22}
\al_k = \gamma_0 \gamma_1^k
\begin{cases}
\abs{\al_1}, & k \text{ odd}, \\ \abs{\al_2}, & k \text{
even},
\end{cases}
\end{align}
where
\begin{equation}
\{\ga_0 = \al_0/\abs{\al_0},\,\ga_1 =
-e^{i(\te_0+\te_1)/2}\}\subset\dD \, \text{ and } \,
\abs{\al_1}\abs{\al_2} = \cos^2((\te_1-\te_0)/4).
\end{equation}
Thus, it remains to show that $\abs{\al_1} = \abs{\al_2}$.
We assume the contrary, $\abs{\al_1} \neq \abs{\al_2}$ and
consider the sequence $\abs{\al} =
\{\abs{\al_k}\}_{k\in\Z}$. Then $\abs{\al}$ is a sequence
of period $2$ Verblunsky coefficients and by \eqref{4.15}
the associated Floquet discriminant, denoted by
$\Delta(\cdot; |\alpha|)$, is given by
\begin{align} \label{5.23}
\De(e^{i\te};\abs{\al}) &= \frac{1}{\rho_1\rho_2}
\left[\frac{e^{i\te}+e^{-i\te}}{2} +
\Re\left(\abs{\al_1}\ol{\abs{\al_2}}\right)\right]
\no \\ &=
\frac{1}{\sqrt{1-\abs{\al_1}^2}\sqrt{1-\abs{\al_2}^2}}
\left[\cos(\te)+\abs{\al_1}{\abs{\al_2}}\right].
\end{align}
Since
\begin{equation}
\si(U_{\abs{\al}}) = \big\{e^{i\te}\in\dD \,\big|\, -1 \leq
\De\big(e^{i\te};\abs{\al}\big) \leq 1 \big\} = \big\{e^{i\te} \in \dD 
\,\big|\,
\la_- \leq \cos(\te) \leq \la_+ \big\},
\end{equation}
where
\begin{align}
\la_\pm = -\abs{\al_1}\abs{\al_2} \pm
\sqrt{1-\abs{\al_1}^2}\sqrt{1-\abs{\al_2}^2}, 
\end{align}
and $\abs{\al_1} \neq \abs{\al_2}$ is equivalent to
$\abs{\la_\pm}<1$, $\si(U_{\abs{\al}})$ consists of two 
arcs. Taking into account
\eqref{5.22}, it follows from Corollary \ref{c3.2} that $\si(U_\al)$ should
also contain two arcs which contradicts the basic hypothesis of 
Theorem \ref{t5.3}. Thus, $\abs{\al_1} = \abs{\al_2}$ and \eqref{5.22}
implies \eqref{5.14}.
\end{proof}

\begin{remark}  \lb{r5.4}
By the last part of Corollary \ref{c3.2}, the phase of $\alpha_0$ in
\eqref{5.3} is a unitary invariant and hence necessarily remains open.
\end{remark}

\appendix
\section{Basic Facts on Caratheodory and Schur Functions}
\lb{A}
\renewcommand{\theequation}{A.\arabic{equation}}
\renewcommand{\thetheorem}{A.\arabic{theorem}}
\setcounter{theorem}{0}
\setcounter{equation}{0}

In this appendix we summarize a few basic properties of Caratheodory
and Schur functions used throughout this manuscript.

We denote by $\D$ and $\dD$ the {\it open unit disk} and the {\it unit
circle} in the complex plane $\C$,
\begin{equation}
\D = \{ z\in\C \st \abs{z} < 1 \}, \quad \dD = \{ \ze\in\C
\st \abs{\ze} = 1 \}
\end{equation}
and by
\begin{equation}
\Cl = \{z\in\C \st \Re(z) < 0\}, \quad
\Cr = \{z\in\C \st \Re(z) > 0\}
\end{equation}
the open left and right complex half-planes, respectively. Moreover,
we orient the unit circle, $\dD$, counterclockwise. By a {\it Laurent
polynomial} we denote a finite linear combination of terms $z^k$,
$k\in\bbZ$, with complex-valued coefficients.

\begin{definition} \lb{dA.1}
Let $f_\pm$, $\varphi_+$, and $1/\varphi_-$ be analytic in $\D$. \\
$(i)$ $f_+$ is called a {\it Caratheodory function} if $f_+\colon
\D\to\Cr$ and $f_-$ is called an {\it anti-Caratheodory function} if
$-f_-$ is a Caratheodory function. \\
$(ii)$ $\varphi_+$ is called a {\it Schur function} if
$\varphi_+\colon\D\to\D$. $\varphi_-$ is called an {\it anti-Schur
function} if $1/\varphi_-$ is a Schur function.
\end{definition}

If $f$ and $g$ are Caratheodory functions, so are
\begin{equation}
f+g, \quad 1/f, \quad \text{and } \, -i\ln(if).
\end{equation}

\begin{theorem} [\cite{Ak65}, Sect.\ 3.1; \cite{AG81}, Sect.\ 69;
\cite{Si04}, Sect.\ 1.3]
${}$ \\
\lb{tA.2} Let $f$ be a Caratheodory function. Then $f$ admits the
Herglotz representation
\begin{align}
& f(z)=ic+ \oint_{\dD} d\mu(\zeta) \, \f{\zeta+z}{\zeta-z},
\quad z\in\D, \lb{A.3} \\
& c=\Im(f(0)), \quad \oint_{\dD} d\mu(\zeta) = \Re(f(0)) < \infty,
\lb{A.4}
\end{align}
where $d\mu$ denotes a nonnegative measure on $\dD$. The measure
$d\mu$ can be reconstructed from $f$ by the formula
\begin{equation}
\mu\big(\Arc\big(\big(e^{i\theta_1},e^{i\theta_2}\big]\big)\big)
=\lim_{\delta\downarrow 0}
\lim_{r\uparrow 1} \f{1}{2\pi}
\oint_{\theta_1+\delta}^{\theta_2+\delta} d\theta \,
\Re\big(f\big(re^{i\theta}\big)\big), \lb{A.6}
\end{equation}
where
\begin{equation} \label{A.5}
\Arc\big(\big(e^{i\theta_1},e^{i\theta_2}\big]\big)
=\big\{e^{i\theta}\in\dD\,|\,
\theta_1<\theta\leq \theta_2\big\}, \quad \theta_1 \in
[0,2\pi), \; \theta_1<\theta_2\leq \theta_1+2\pi.
\end{equation}
Conversely, the right-hand side of \eqref{A.3} with $c\in\bbR$ and
$d\mu$ a finite $($nonnegative$)$ measure on $\dD$ defines a
Caratheodory function.
\end{theorem}

We note that additive nonnegative constants on the
right-hand side of \eqref{A.3} can be absorbed into the measure
$d\mu$ since
\begin{equation}
\oint_\dD d\mu_0(\zeta) \, \f{\zeta+z}{\zeta-z}=1, \lb{A.5a}
\end{equation}
where
\begin{equation}
d\mu_0(\zeta)=\f{d\theta}{2\pi}, \quad \zeta=e^{i\theta}, \;
\theta\in [0,2\pi] \lb{A.5b}
\end{equation}
denotes the normalized Lebesgue measure on the unit circle $\dD$.

A useful fact on Caratheodory functions $f$ is a certain
monotonicity property they exhibit on open connected arcs
of the unit circle away from the support of the measure
$d\mu$ in the Herglotz representation \eqref{A.3}. More
precisely, suppose
$\Arc\big(\big(e^{i\theta_1},e^{i\theta_2}\big)\big)\subset
(\dD\backslash\supp(d\mu))$, $\theta_1<\theta_2$, then $f$
has an analytic continuation through
$\Arc\big(\big(e^{i\theta_1},e^{i\theta_2}\big)\big)$ and
it is purely imaginary on
$\Arc\big(\big(e^{i\theta_1},e^{i\theta_2}\big)\big)$.
Moreover,
\begin{equation}
\f{d}{d\theta}f\big(e^{i\theta}\big)=-\f{i}{2}
\int_{[0,2\pi]\backslash(\theta_1,\theta_2)}
d\mu\big(e^{it}\big) \f{1}{\sin^2((t-\theta)/2)}, \quad
\theta\in(\theta_1,\theta_2).  \lb{A.9a}
\end{equation}
In particular,
\begin{equation}
-i\f{d}{d\theta}f\big(e^{i\theta}\big)<0, \quad
\theta\in(\theta_1,\theta_2).
\end{equation}

We recall that any Caratheodory function $f$ has finite radial limits to
the unit circle $\mu_0$-almost everywhere, that is,
\begin{equation}
f(\zeta)=\lim_{r\uparrow 1} f(r\zeta) \, \text{ exists and is finite for
$\mu_0$-a.e.\ $\zeta\in\dD$.}
\end{equation}

Next, we consider the decomposition of the measure $d\mu$ in the
Herglotz representation \eqref{A.3} of the Caratheodory
function $f$ into its absolutely continuous and singular parts,
$d\mu_{\rm ac}$ and $d\mu_{\rm s}$, respectively,
\begin{equation}
d\mu=d\mu_{\rm ac}+d\mu_{\rm s}.  \lb{A.}
\end{equation}
The absolutely continuous part $d\mu_{\rm ac}$ of $d\mu$ is given by
\begin{equation}
d\mu_{\rm ac}(\zeta) = \lim_{r\uparrow 1}
\Re(f(r\zeta)) \, d\mu_0(\zeta), \quad \zeta \in\dD.
\lb{A.10}
\end{equation}
The set
\begin{equation}
S_{\mu_{\rm ac}}=\{\zeta\in\dD\,|\, \lim_{r\uparrow
1}\Re(f(r\zeta)) =\Re(f(\zeta))>0 \text{ exists finitely}\}
\lb{A.10a}
\end{equation}
is an essential support of $d\mu_{\rm ac}$ and its essential closure,
$\ol{S_{\mu_{\rm ac}}}^e$, coincides with the topological support,
$\supp(d\mu_{\rm ac})$ (the smallest closed support), of
$d\mu_{\rm ac}$,
\begin{equation}
\ol{S_{\mu_{\rm ac}}}^e=\supp \, (d\mu_{\rm ac}).  \lb{A.10b}
\end{equation}
Moreover, the set
\begin{align}
S_{\mu_{\rm s}} = \{\zeta\in\dD\,|\, \lim_{r\uparrow
1}\Re(f(r\zeta)) = \infty \} \label{A.10d}
\end{align}
is an essential support of the singular part $d\mu_{\rm s}$
of the measure $d\mu$, and
\begin{equation}
\lim_{r\uparrow 1} (1-r)f(r\zeta)=\lim_{r\uparrow 1}
(1-r)\Re(f(r\zeta))\geq 0 \, \text{ exists for all $\zeta\in\dD$.}
\lb{A.10C}
\end{equation}
In particular, $\zeta_0\in\dD$ is a pure point of $d\mu$ if and only if
\begin{equation}
\mu(\{\zeta_0\})=\lim_{r \uparrow 1} \left(\f{1-r}{2}\right)
f(r\zeta_0)>0. \lb{A.10c}
\end{equation}

Given a Caratheodory (resp., anti-Caratheodory) function $f_+$
(resp. $f_-$) defined in $\D$ as in \eqref{A.3}, one extends $f_\pm$ to
all of $\bbC\backslash\dD$ by
\begin{equation}
f_\pm(z)=ic_\pm \pm \oint_{\dD} d\mu_\pm (\zeta) \,
\f{\zeta+z}{\zeta-z}, \quad z\in\bbC\backslash\dD, \;
c_\pm\in\bbR. \lb{A.6A}
\end{equation}
In particular,
\begin{equation}
f_\pm(z) = -\ol{f_\pm(1/\ol{z})}, \quad
z\in\C\backslash\ol{\D}. \lb{A.7}
\end{equation}
Of course, this continuation of $f_\pm|_{\D}$ to
$\bbC\backslash\ol\D$, in general, is not an analytic
continuation of $f_\pm|_\D$. With $f_\pm$ defined on
$\bbC\backslash\dD$ by \eqref{A.6} one infers the mapping
properties
\begin{equation}
f_+\colon \D \to \Cr, \quad f_+\colon\bbC\backslash
\ol\D\to\Cl, \quad f_-\colon \D \to \Cl, \quad f_-\colon\bbC\backslash
\ol\D\to\Cr.
\end{equation}

Next, given the functions $f_\pm$ defined in $\bbC\backslash\dD$
as in \eqref{A.6}, we introduce the functions $\varphi_\pm$
by
\begin{equation}
\varphi_\pm(z)=\f{f_\pm(z)-1}{f_\pm(z)+1}, \quad
z\in\bbC\backslash\dD.  \lb{A.18a}
\end{equation}
Then $\varphi_\pm$ have the mapping properties
\begin{align}
\begin{split}
& \varphi_+\colon\D\to\D, \quad
1/\varphi_+\colon\bbC\backslash\ol\D \to\D \quad
(\varphi_+\colon\bbC\backslash\ol\D\to
(\bbC\backslash\ol\D)\cup\{\infty\}), \\ &
\varphi_-\colon\bbC\backslash\ol\D\to\D, \quad
1/\varphi_-\colon\D\to\D \quad (\varphi_-\colon\D\to
(\bbC\backslash\ol\D)\cup\{\infty\}),
\end{split}
\end{align}
in particular, $\varphi_+|_{\D}$ (resp., $\varphi_-|_{\D}$)
is a Schur (resp., anti-Schur) function. Moreover,
\begin{equation}
f_\pm(z)=\f{1+\varphi_\pm (z)}{1-\varphi_\pm (z)}, \quad
z\in\bbC\backslash\dD.
\end{equation}

In analogy to the exponential representation of Nevanlinna--Herglotz
functions (i.e., functions analytic in the open complex upper
half-plane $\bbC_+$ with a strictly positive imaginary part on
$\bbC_+$, cf.\ \cite{AD56}, \cite{AD64}, \cite{GT00}, \cite{KK74}) one
obtains the following result.

\begin{theorem} \lb{tA.3}
Let $f$ be a Caratheodory function. Then $-i\ln(if)$ is a
Caratheodory function and $f$ has the exponential Herglotz
representation,
\begin{align}
& -i\ln(if(z))=id+ \oint_{\dD} d\mu_0(\zeta) \,\Upsilon (\zeta) \,
\f{\zeta+z}{\zeta-z},
\quad z\in\D, \lb{A.11} \\
& \;\, d=-\Re(\ln(f(0))), \quad 0 \leq \Upsilon (\zeta)\leq \pi
\, \text{ for $\mu_0$-a.e.\ $\zeta\in\dD$}.
\lb{A.13}
\end{align}
$\Upsilon$ can be reconstructed from $f$ by
\begin{align}
\begin{split}
\Upsilon (\zeta) &= \lim_{r\uparrow 1}\Re[-i\ln(if(r\zeta))]  \\
& =(\pi/2)+\lim_{r\uparrow 1}\Im[\ln(f(r\zeta))]\, \text{ for
$\mu_0$-a.e.\ $\zeta\in\dD$.}
\end{split}
\end{align}
\end{theorem}

Next we briefly turn to matrix-valued Caratheodory
functions. We denote as usual $\Re(A)=(A+A^*)/2$,
$\Im(A)=(A-A^*)/(2i)$, etc., for square matrices $A$.

\begin{definition} \lb{dA.4}
Let $m\in\bbN$ and $\cF$ be an $m\times m$ matrix-valued function
analytic in $\D$. $\cF$ is called a {\it Caratheodory matrix} if
$\Re(\cF(z))\geq 0$ for all $z\in\D$.
\end{definition}

\begin{theorem} \label{tA.5}
Let $\cF$ be an $m\times m$ Caratheodory matrix, $m\in\bbN$.
Then $\cF$ admits the Herglotz representation
\begin{align}
& \cF(z)=iC+ \oint_{\dD} d\Omega(\zeta) \,
\f{\zeta+z}{\zeta-z}, \quad z\in\D, \lb{A.24a}
\\
& C=\Im(\cF(0)), \quad \oint_{\dD} d\Omega(\zeta) =
\Re(\cF(0)), \lb{A.25a}
\end{align}
where $d\Omega$ denotes a nonnegative $m \times m$
matrix-valued measure on $\dD$. The measure $d\Omega$ can
be reconstructed from $\cF$ by the formula
\begin{align}
&\Omega\big(\Arc\big(\big(e^{i\theta_1},e^{i\theta_2}\big]\big)\big)
=\lim_{\delta\downarrow 0} \lim_{r\uparrow 1} \f{1}{2\pi}
\oint_{\theta_1+\delta}^{\theta_2+\delta} d\theta \,
\Re\big(\cF\big(re^{i\theta}\big)\big), \\
& \hspace*{3.9cm} \theta_1 \in
[0,2\pi), \; \theta_1<\theta_2\leq \theta_1+2\pi. \no
\end{align}
Conversely, the right-hand side of equation \eqref{A.24a}
with $C = C^*$ and $d\Omega$ a finite nonnegative $m \times
m$ matrix-valued measure on $\dD$ defines a Caratheodory
matrix.
\end{theorem}

\section{A Summary of Weyl--Titchmarsh Theory for CMV Operators on
Half-Lattices and on $\bbZ$}
\lb{B}
\renewcommand{\theequation}{B.\arabic{equation}}
\renewcommand{\thetheorem}{B.\arabic{theorem}}
\setcounter{theorem}{0}
\setcounter{equation}{0}

We start by introducing some of the basic notations used throughout
this paper. Detailed proofs of all facts in this appendix (and a lot of
additional material) can be found in \cite{GZ05}.

In the following, let $\ltz$ be the usual Hilbert space of all square
summable complex-valued sequences with scalar product
$(\cdot,\cdot)_{\ell^2(\bbZ)}$ linear in the second argument. The
{\it standard basis} in $\ell^2(\bbZ)$ is denoted by
\begin{equation}
\{\delta_k\}_{k\in\bbZ}, \quad
\delta_k=(\dots,0,\dots,0,\underbrace{1}_{k},0,\dots,0,\dots)^\top,
\; k\in\bbZ.
\end{equation}
For $J\subseteq\bbR$ an interval, we will identify
$\ell^2(J\cap\bbZ)\oplus\ell^2(J\cap\bbZ)$ and
$\ell^2(J\cap\bbZ)\otimes\bbC^2$ and then use the simplified notation
$\ell^2(J\cap\bbZ)^2$. For simplicity, the identity operator on
$\ell^2(J\cap\bbZ)$ is abbreviated by $I$ without separately
indicating its dependence on $J$.

Throughout this appendix we make the following basic assumption:

\begin{hypothesis} \lb{hA.4}
Let $\alpha$ be a sequence of complex numbers such that
\begin{equation} \label{A.15}
\alpha=\{\al_k\}_{k \in \Z} \subset \D.
\end{equation}
\end{hypothesis}

Given a sequence $\alpha$ satisfying \eqref{A.15}, we define
the sequence of positive real numbers $\{\rho_k\}_{k\in\bbZ}$ and two
sequences of complex numbers with positive real parts
$\{a_k\}_{k\in\bbZ}$ and $\{b_k\}_{k\in\bbZ}$ by
\begin{equation}
\rho_k = \sqrt{1-\abs{\al_k}^2}, \quad a_k = 1+\al_k, \quad
b_k = 1-\al_k, \quad k \in \Z.
\end{equation}

Following Simon \cite{Si04}, we call $\alpha_k$ the Verblunsky
coefficients in honor of Verblunsky's pioneering work in the theory
of orthogonal polynomials on the unit circle \cite{Ve35}, \cite{Ve36}.

Next, we also introduce a sequence of $2\times 2$ unitary
matrices $\te_k$ by
\begin{equation} \label{A.19}
\te_k = \begin{pmatrix} -\al_k & \rho_k \\ \rho_k &
\ol{\al_k}
\end{pmatrix},
\quad k \in \Z,
\end{equation}
and two unitary operators $V$ and $W$ on $\ltz$ by their
matrix representations in the standard basis of $\ell^2(\bbZ)$ as
follows,
\begin{align} \label{A.20}
V &= \begin{pmatrix} \ddots & & &
\raisebox{-3mm}[0mm][0mm]{\hspace*{-5mm}\Huge $0$}  \\ &
\te_{2k-2} & & \\ & & \te_{2k} & & \\ &
\raisebox{0mm}[0mm][0mm]{\hspace*{-10mm}\Huge $0$} & &
\ddots
\end{pmatrix}, \quad
W &= \begin{pmatrix} \ddots & & &
\raisebox{-3mm}[0mm][0mm]{\hspace*{-5mm}\Huge $0$}
\\ & \te_{2k-1} &  &  \\ &  & \te_{2k+1} &  & \\ &
\raisebox{0mm}[0mm][0mm]{\hspace*{-10mm}\Huge $0$} & &
\ddots
\end{pmatrix},
\end{align}
where
\begin{align}
\begin{pmatrix}
V_{2k-1,2k-1} & V_{2k-1,2k} \\
V_{2k,2k-1}   & V_{2k,2k}
\end{pmatrix} =  \te_{2k},
\quad
\begin{pmatrix}
W_{2k,2k} & W_{2k,2k+1} \\ W_{2k+1,2k}  & W_{2k+1,2k+1}
\end{pmatrix} =  \te_{2k+1},
\quad k\in\Z.
\end{align}
Moreover, we introduce the unitary operator $U$ on $\ltz$ by
\begin{equation} \label{A.22}
U = VW,
\end{equation}
or in matrix form, in the standard basis of $\ell^2(\bbZ)$, by
\begin{align}
U = \begin{pmatrix} \ddots &&\hspace*{-8mm}\ddots
&\hspace*{-10mm}\ddots &\hspace*{-12mm}\ddots
&\hspace*{-14mm}\ddots &&&
\raisebox{-3mm}[0mm][0mm]{\hspace*{-6mm}{\Huge $0$}}
\\
&0& -\al_{0}\rho_{-1} & -\ol{\al_{-1}}\al_{0} & -\al_{1}\rho_{0} &
\rho_{0}\rho_{1}
\\
&& \rho_{-1}\rho_{0} &\ol{\al_{-1}}\rho_{0} &
-\ol{\al_{0}}\al_{1} & \ol{\al_{0}}\rho_{1} & 0
\\
&&0& -\al_{2}\rho_{1} & -\ol{\al_{1}}\al_{2} &
-\al_{3}\rho_{2} & \rho_{2}\rho_{3}
\\
&\raisebox{-4mm}[0mm][0mm]{\hspace*{-6mm}{\Huge $0$}} &&
\rho_{1}\rho_{2} & \ol{\al_{1}}\rho_{2} & -\ol{\al_{2}}\al_{3} &
\ol{\al_{2}}\rho_{3} &0&
\\
&&&&\hspace*{-14mm}\ddots &\hspace*{-14mm}\ddots
&\hspace*{-14mm}\ddots &\hspace*{-8mm}\ddots &\ddots
\end{pmatrix}. \lb{A.23}
\end{align}
Here terms of the form $-\ol{\alpha_k}\alpha_{k+1}$, $k\in\Z$,
represent the diagonal entries in the infinite matrix \eqref{A.23}.
We will call the operator $U$ on $\ell^2(\bbZ)$ the CMV operator since
\eqref{A.19}--\eqref{A.23} in the context of the semi-infinite (i.e.,
half-lattice) case were first obtained by Cantero, Moral, and
Vel\'azquez in \cite{CMV03}.

Finally, let $\UU$ denote the unitary operator on
$\ltz^2$ defined by
\begin{equation} \label{A.24}
\UU = \begin{pmatrix}
U & 0 \\
0 & U^\top
\end{pmatrix}
=
\begin{pmatrix}
VW & 0 \\
0 & WV
\end{pmatrix}
=
\begin{pmatrix}
0 & V \\
W & 0
\end{pmatrix}^2.
\end{equation}
One observes remnants of a certain ``supersymmetric'' structure
in
$\left(\begin{smallmatrix} 0 & V \\ W & 0 \end{smallmatrix}\right)$
which is also reflected in the following result.

\begin{lemma} \label{lA.1}
Let $z\in\bbC\backslash\{0\}$. Then the following items $(i)$--$(iv)$ are
equivalent:
\begin{align}
& (i) \quad U u(z,\cdot) = z u(z,\cdot). \\
& (ii) \quad \UU \binom{u(z,\cdot)}{v(z,\cdot)} =
z \binom{u(z,\cdot)}{v(z,\cdot)}, \quad
(W u)(z,\cdot)=z v(z,\cdot). \\
& (iii) \quad \UU \binom{u(z,\cdot)}{v(z,\cdot)} =
z \binom{u(z,\cdot)}{v(z,\cdot)}, \quad
(V v)(z,\cdot)=u(z,\cdot). \\
& (iv) \quad \binom{u(z,k)}{v(z,k)} = T(z,k) \binom{u(z,k-1)}{v(z,k-1)},
\quad k\in\Z,  \label{A.28}
\end{align}
where the transfer matrices $T(z,k)$, $z\in\Cz$, $k\in\Z$,
are given by
\begin{equation}
T(z,k) = \begin{cases} \frac{1}{\rho_{k}} \begin{pmatrix}
\al_{k} & z \\ 1/z & \ol{\al_{k}} \end{pmatrix},  &
\text{$k$ odd,}  \\ \frac{1}{\rho_{k}} \begin{pmatrix}
\ol{\al_{k}} & 1 \\ 1 & \al_{k} \end{pmatrix}, & \text{$k$
even.} \end{cases}     \label{A.31}
\end{equation}
\end{lemma}

If one sets $\al_{k_0} = e^{is}$, $s\in [0,2\pi)$, for some
reference point $k_0\in\Z$, then the operator $U$ splits
into a direct sum of two half-lattice operators
$U_{-,k_0-1}^{(s)}$ and $U_{+,k_0}^{(s)}$ acting on
$\ell^2((-\infty,k_0-1]\cap\Z)$ and on
$\ell^2([k_0,\infty)\cap\Z)$, respectively. Explicitly, one
obtains
\begin{align}
\begin{split}
& U=U_{-,k_0-1}^{(s)} \oplus U_{+,k_0}^{(s)} \, \text{ in }
\, \ell^2((-\infty,k_0-1]\cap\Z) \oplus
\ell^2([k_0,\infty)\cap\Z) \\ & \text{if } \, \al_{k_0} =
e^{is}, \; s\in [0,2\pi). \lb{B.18}
\end{split}
\end{align}
Similarly, one obtains $W_{-,k_0-1}^{(s)}$,
$V_{-,k_0-1}^{(s)}$ and $W_{+,k_0}^{(s)}$,
$V_{+,k_0}^{(s)}$ such that
\begin{equation}
U_{\pm,k_0}^{(s)} = V_{\pm,k_0}^{(s)} W_{\pm,k_0}^{(s)}.
\end{equation}
For simplicity we will abbreviate
\begin{equation}
U_{\pm,k_0} =
U_{\pm,k_0}^{(s=0)}=V_{\pm,k_0}^{(s=0)}W_{\pm,k_0}^{(s=0)}
=V_{\pm,k_0} W_{\pm,k_0}.  \lb{B.17}
\end{equation}
In addition, we introduce on
$\ell^2([k_0,\pm\infty)\cap\Z)^2$ the half-lattice
operators $\UU_{\pm,k_0}^{(s)}$ by
\begin{equation}
\UU_{\pm,k_0}^{(s)} = \begin{pmatrix} U_{\pm,k_0}^{(s)} & 0
\\ 0 & (U_{\pm,k_0}^{(s)})^\top
\end{pmatrix}
=\begin{pmatrix} V_{\pm,k_0}^{(s)} W_{\pm,k_0}^{(s)} & 0 \\
0 & W_{\pm,k_0}^{(s)} V_{\pm,k_0}^{(s)}
\end{pmatrix}.
\end{equation}
By $\UU_{\pm,k_0}$ we denote the half-lattice operators
defined for $s=0$,
\begin{equation}
\UU_{\pm,k_0} = \UU_{\pm,k_0}^{(s=0)} = \begin{pmatrix}
U_{\pm,k_0} & 0 \\ 0 & (U_{\pm,k_0})^\top
\end{pmatrix}
=
\begin{pmatrix}
V_{\pm,k_0}W_{\pm,k_0} & 0 \\
0 & W_{\pm,k_0}V_{\pm,k_0}
\end{pmatrix}.
\end{equation}

\begin{lemma} \label{lA.2}
Let $z\in\bbC\backslash\{0\}$ and $k_0\in\bbZ$. Consider
sequences $\{\hatt p_+(z,k,k_0)\}_{k\geq k_0}$, $\{\hatt
r_+(z,k,k_0)\}_{k\geq k_0}$. Then, the following items
$(i)$--$(iii)$ are equivalent:
\begin{align}
& (i) \quad \UU_{+,k_0} \binom{\hatt
p_+(z,\cdot,k_0)}{\hatt r_+(z,\cdot,k_0)} = z \binom{\hatt
p_+(z,\cdot,k_0)}{\hatt r_+(z,\cdot,k_0)}, \quad
W_{+,k_0}\hatt p_+(z,\cdot,k_0) = z \hatt r_+(z,\cdot,k_0).
\\[1mm] & (ii) \quad \UU_{+,k_0} \binom{\hatt
p_+(z,\cdot,k_0)}{\hatt r_+(z,\cdot,k_0)} = z \binom{\hatt
p_+(z,\cdot,k_0)}{\hatt r_+(z,\cdot,k_0)}, \quad
V_{+,k_0}\hatt r_+(z,\cdot,k_0) = \hatt p_+(z,\cdot,k_0).
\\ & (iii) \quad \binom{\hatt p_+(z,k,k_0)}{\hatt
r_+(z,k,k_0)} = T(z,k) \binom{\hatt p_+(z,k-1,k_0)}{\hatt
r_+(z,k-1,k_0)}, \quad k > k_0, \label{A.41}  \\
& \phantom{(vi)} \quad \;\; \text{ assuming } \, 
\hatt p_+(z,k_0,k_0) = \begin{cases}
z\hatt r_+(z,k_0,k_0), & \text{$k_0$ odd}, \\ \hatt
r_+(z,k_0,k_0), & \text{$k_0$ even}. \end{cases} \lb{A.42}
\end{align}
Next, consider sequences $\{\hatt p_-(z,k,k_0)\}_{k\leq
k_0}$, $\{\hatt r_-(z,k,k_0)\}_{k\leq k_0}$. Then, the
following items $(iv)$--$(vi)$are equivalent:
\begin{align}
& (iv) \quad \UU_{-,k_0} \binom{\hatt
p_-(z,\cdot,k_0)}{\hatt r_-(z,\cdot,k_0)} = z \binom{\hatt
p_-(z,\cdot,k_0)}{\hatt r_-(z,\cdot,k_0)}, \quad
W_{-,k_0}\hatt p_-(z,\cdot,k_0) = z \hatt r_-(z,\cdot,k_0).
\\[1mm] & (v) \quad  \UU_{-,k_0} \binom{\hatt
p_-(z,\cdot,k_0)}{\hatt r_-(z,\cdot,k_0)} = z \binom{\hatt
p_-(z,\cdot,k_0)}{\hatt r_-(z,\cdot,k_0)}, \quad
V_{-,k_0}\hatt r_-(z,\cdot,k_0) = \hatt p_-(z,\cdot,k_0).
\\ & (vi) \quad \binom{\hatt p_-(z,k-1),k_0}{\hatt
r_-(z,k-1,k_0)}=T(z,k)^{-1} \binom{\hatt
p_-(z,k,k_0)}{\hatt r_-(z,k,k_0)}, \quad k \leq k_0,
\label{A.46} \\ 
& \phantom{(xii)} \quad  \text{ assuming } \, 
\hatt p_-(z,k_0,k_0) =\begin{cases} -\hatt r_-(z,k_0,k_0),
& \text{$k_0$ odd,}
\\ -z\hatt r_-(z,k_0,k_0), & \text{$k_0$ even.} \end{cases}
\lb{A.47}
\end{align}
\end{lemma}

In the following, we denote by
$\Big(\begin{smallmatrix}p_+(z,k,k_0)\\
r_+(z,k,k_0)\end{smallmatrix}\Big)_{k\geq k_0}$ and
$\Big(\begin{smallmatrix}q_+(z,k,k_0)\\
s_+(z,k,k_0)\end{smallmatrix}\Big)_{k\geq k_0}$,
$z\in\bbC\backslash\{0\}$, two linearly independent solutions of
\eqref{A.41} with the following initial conditions:
\begin{align}
\binom{p_+(z,k_0,k_0)}{r_+(z,k_0,k_0)} = \begin{cases} \binom{z}{1}, &
\text{$k_0$ odd,} \\[1mm]
\binom{1}{1}, & \text{$k_0$ even,} \end{cases} \quad
\binom{q_+(z,k_0,k_0)}{s_+(z,k_0,k_0)} = \begin{cases} \binom{z}{-1}, &
\text{$k_0$ odd,} \\[1mm]
\binom{-1}{1}, & \text{$k_0$ even.} \end{cases} \lb{A.50}
\end{align}
Similarly, we denote by
$\Big(\begin{smallmatrix}p_-(z,k,k_0)\\
r_-(z,k,k_0)\end{smallmatrix}\Big)_{k\leq k_0}$ and
$\Big(\begin{smallmatrix}q_-(z,k,k_0)\\
s_-(z,k,k_0)\end{smallmatrix}\Big)_{k\leq k_0}$,
$z\in\bbC\backslash\{0\}$, two linearly independent solutions of
\eqref{A.46} with the following initial conditions:
\begin{align}
\binom{p_-(z,k_0,k_0)}{r_-(z,k_0,k_0)} = \begin{cases} \binom{1}{-1}, &
\text{$k_0$ odd,} \\[1mm]
\binom{-z}{1}, & \text{$k_0$ even,} \end{cases} \quad
\binom{q_-(z,k_0,k_0)}{s_-(z,k_0,k_0)} = \begin{cases} \binom{1}{1}, &
\text{$k_0$ odd,} \\[1mm]
\binom{z}{1}, & \text{$k_0$ even.} \end{cases} \lb{A.52}
\end{align}
Using \eqref{A.28} one extends
$\Big(\begin{smallmatrix}p_+(z,k,k_0)\\r_+(z,k,k_0)
\end{smallmatrix}\Big)_{k\geq k_0}$,
$\Big(\begin{smallmatrix}q_+(z,k,k_0)\\s_+(z,k,k_0)
\end{smallmatrix}\Big)_{k\geq k_0}$, $z\in\bbC\backslash\{0\}$,
to $k < k_0$. In the same manner, one extends
$\Big(\begin{smallmatrix}p_-(z,k,k_0)\\r_-(z,k,k_0)
\end{smallmatrix}\Big)_{k\leq k_0}$ and
$\Big(\begin{smallmatrix}q_-(z,k,k_0)\\s_-(z,k,k_0)
\end{smallmatrix}\Big)_{k\leq k_0}$, $z\in\bbC\backslash\{0\}$,
to $k> k_0$. These extensions will be denoted by
$\Big(\begin{smallmatrix}p_\pm(z,k,k_0)\\r_\pm
(z,k,k_0)\end{smallmatrix}\Big)_{k\in\Z}$ and
$\Big(\begin{smallmatrix}q_\pm (z,k,k_0)\\s_\pm
(z,k,k_0)\end{smallmatrix}\Big)_{k\in\Z}$.\ Moreover, it
follows from \eqref{A.28} that $p_\pm(z,k,k_0)$,
$q_\pm(z,k,k_0)$, $r_\pm(z,k,k_0)$, and $s_\pm(z,k,k_0)$,
$k,k_0\in\Z$, are Laurent polynomials in $z$.

\begin{remark} \label{rA.3}
By Lemma \ref{lA.2},
$\Big(\begin{smallmatrix}p_\pm(z,k,k_0)\\
r_\pm(z,k,k_0)\end{smallmatrix}\Big)_{k\gtreqless k_0}$,
$z\in\bbC\backslash\{0\}$, $k_0\in\Z$, are generalized
eigenvectors of the operators $\UU_{\pm,k_0}$. Moreover,
by Lemma \ref{lA.1}, $\Big(\begin{smallmatrix}p_\pm
(z,k,k_0)\\r_\pm (z,k,k_0)\end{smallmatrix}\Big)_{k\in\Z}$
and $\Big(\begin{smallmatrix}q_\pm (z,k,k_0)\\s_\pm
(z,k,k_0)\end{smallmatrix}\Big)_{k\in\Z}$,
$z\in\bbC\backslash\{0\}$, $k_0\in\Z$, are generalized
eigenvectors of $\UU$.
\end{remark}

\begin{lemma}
Let $k_0\in\bbZ$. Then the sets of Laurent polynomials
$\{p_\pm(\cdot,k,k_0)\}_{k\gtreqless k_0}$ and
$\{r_\pm(\cdot,k,k_0)\}_{k\gtreqless k_0}$ form complete
orthonormal systems in $\Lt{\pm}$, where
\begin{equation}
d\mu_\pm(\zeta,k_0) = d(\de_{k_0},E_{U_{\pm,k_0}}(\zeta)\de_{k_0}),
\quad \zeta\in\dD, \label{A.53}
\end{equation}
and $dE_{U_{\pm,k_0}}(\cdot)$ denote the operator-valued
spectral measures of the operators $U_{\pm,k_0}$,
\begin{equation}
U_{\pm,k_0}=\oint_{\dD} dE_{U_{\pm,k_0}}(\zeta)\,\zeta.
\end{equation}
\end{lemma}

We note that the measures $d\mu_\pm(\cdot,k_0)$, $k_0\in\bbZ$, are
not only nonnegative but also supported on an infinite set.

\begin{remark} \label{rB.6}
In connection with our introductory remarks in \eqref{1.3}--\eqref{1.4e}
we note that $d\sigma_+=d\mu_+(\cdot,0)$ and that in accordance with
\cite[Proposition 4.2.2]{Si04}
\begin{align}
\begin{split}
p_+(\zeta,k,0)&=\begin{cases} \gamma_k \zeta^{-(k-1)/2}
\varphi_+(\zeta,k), & \text{$k$ odd,} \\[1mm]
\gamma_k \zeta^{-k/2} \varphi^*_+(\zeta,k), & \text{$k$ even,}
\end{cases}  \\
r_+(\zeta,k,0)&=\begin{cases} \gamma_{k} \zeta^{-(k+1)/2}
\varphi^*_+(\zeta,k), & \text{$k$ odd,} \\[1mm]
\gamma_k \zeta^{-k/2} \varphi_+(\zeta,k), & \text{$k$ even;}
\end{cases} \quad \zeta\in\dD.
\end{split}
\end{align}
\end{remark}

\begin{corollary} \label{cA.4}
Let $k_0\in\bbZ$. Then the half-lattice CMV operators
$U_{\pm,k_0}$ are unitarily equivalent to the operators of
multiplication by the function $id$ $($where $id(\ze)=\ze$,
$\ze\in\dD$$)$ on $\Lt{\pm}$. In particular,
\begin{equation}
\si(U_{\pm,k_0}) =\supp \, (d\mu_\pm(\cdot,k_0))
\end{equation}
and the spectrum of $U_{\pm,k_0}$ is simple.
\end{corollary}

\begin{lemma} \lb{cA.10}
Let $k_0\in\bbZ$. Then the sets of two-dimensional Laurent polynomials
$\Big(\begin{smallmatrix}p_\pm(z,k,k_0)\\
r_\pm(z,k,k_0)\end{smallmatrix}\Big)_{k\gtreqless k_0}$ and
$\Big(\begin{smallmatrix}q_\pm(z,k,k_0)\\
s_\pm(z,k,k_0)\end{smallmatrix}\Big)_{k\gtreqless k_0}$
satisfy the relation
\begin{equation}
\binom{q_\pm(z,\cdot,k_0)}{s_\pm(z,\cdot,k_0)} + m_\pm(z,k_0)
\binom{p_\pm(z,\cdot,k_0)}{r_\pm(z,\cdot,k_0)} \in
\ell^2([k_0,\pm\infty)\cap\Z)^2, \quad
z\in\bbC\backslash(\dD\cup\{0\}), \label{A.63}
\end{equation}
where the coefficients $m_\pm(z,k_0)$ are given by
\begin{align}
m_\pm(z,k_0) &= \pm
(\delta_{k_0},(U_{\pm,k_0}+zI)(U_{\pm,k_0}-zI)^{-1}
\delta_{k_0})_{\ell^2(\Z)} \lb{A.66} \\
& =\pm \oint_\dD d\mu_{\pm}(\zeta,k_0)\,
\frac{\zeta+z}{\zeta-z},
\quad z\in\bbC\backslash\dD \lb{A.67}
\intertext{with}
m_\pm(0,k_0)&=\pm\oint_{\dD} d\mu_\pm(\zeta,k_0)=\pm 1.
\end{align}
\end{lemma}

\begin{lemma} \label{lA.8}
Let $k_0\in\bbZ$. Then relation \eqref{A.63} uniquely
determines the functions $m_\pm(\cdot,k_0)$ on
$\bbC\backslash\dD$.
\end{lemma}

\begin{theorem} \lb{tA.10}
Let $k_0\in\bbZ$. Then there exist
unique functions $M_\pm(\cdot,k_0)$ such that
\begin{align}
&\binom{u_\pm(z,\cdot,k_0)}{v_\pm(z,\cdot,k_0)} =
\binom{q_+(z,\cdot,k_0)}{s_+(z,\cdot,k_0)} + M_\pm(z,k_0)
\binom{p_+(z,\cdot,k_0)}{r_+(z,\cdot,k_0)} \in
\ell^2([k_0,\pm\infty)\cap\Z)^2, \no \\
& \hspace*{8cm} z\in\bbC\backslash(\dD\cup\{0\}). \label{A.68}
\end{align}
\end{theorem}

We will call $u_\pm(z,\cdot,k_0)$ (resp.,
$v_\pm(z,\cdot,k_0)$) {\it Weyl--Titchmarsh solutions} of
$U$ (resp., $U^\top$). Similarly, we will call
$m_\pm(z,k_0)$ as well as $M_\pm(z,k_0)$ the {\it
half-lattice Weyl--Titchmarsh $m$-functions} associated
with $U_{\pm,k_0}$. (See also \cite{Si04a} for a comparison
of various alternative notions of Weyl--Titchmarsh
$m$-functions for $U_{+,k_0}$.)

One verifies that
\begin{align}
M_+(z,k_0) &= m_+(z,k_0), \quad z\in\bbC\backslash\dD,
\lb{A.69} \\ M_+(0,k_0) &=1, \lb{A.70} \\ M_-(z,k_0) &=
\frac{\Re(a_{k_0}) +
i\Im(b_{k_0})m_-(z,k_0-1)}{i\Im(a_{k_0}) +
\Re(b_{k_0})m_-(z,k_0-1)}, \quad z\in\bbC\backslash\dD,
\lb{A.71} \\ M_-(0,k_0)
&=\f{\alpha_{k_0}+1}{\alpha_{k_0}-1}. \lb{A.72}
\end{align}
In particular, one infers that $M_\pm$ are analytic at $z=0$.

\begin{lemma}
Let $k\in\bbZ$. Then the functions $M_+(\cdot,k)|_{\D}$ $($resp.,
$M_-(\cdot,k)|_{\D}$$)$ are Caratheodory $($resp.,
anti-Caratheodory\,$)$ functions. Moreover, $M_\pm$ satisfy the
following Riccati-type equation
\begin{align}
&(z\ol{b_k}-b_k)M_\pm(z,k-1)M_\pm(z,k)+(z\ol{b_k}+b_k)M_\pm(z,k)
-(z\ol{a_k}+a_k)M_\pm(z,k-1) \no \\
& \quad =z\ol{a_k}-a_k, \quad z\in\bbC\backslash\dD. \label{A.77}
\end{align}
\end{lemma}

In addition, we introduce the functions $\Phi_\pm(\cdot,k)$,
$k\in\bbZ$, by
\begin{align}
\Phi_\pm(z,k) = \frac{M_\pm(z,k)-1}{M_\pm(z,k)+1},
\quad z\in\C\backslash\dD. \lb{A.78}
\end{align}
One then verifies,
\begin{equation}
M_\pm(z,k) = \frac{1+\Phi_\pm(z,k)}{1-\Phi_\pm(z,k)},
\quad z\in\C\backslash\dD. \lb{A.79}
\end{equation}

\begin{lemma}
Let $z\in\C\backslash(\dD\cup\{0\})$, $k_0, k\in\Z$. Then
the functions $\Phi_\pm(\cdot,k)$ satisfy the following
equalities
\begin{equation}
\Phi_\pm(z,k) = \begin{cases}
z\frac{v_\pm(z,k,k_0)}{u_\pm(z,k,k_0)}, &\text{$k$ odd,}
\\
\frac{u_\pm(z,k,k_0)}{v_\pm(z,k,k_0)}, & \text{$k$ even,}
\end{cases} \lb{B.72}
\end{equation}
where $u_\pm(\cdot,k,k_0)$ and $v_\pm(\cdot,k,k_0)$ are the
polynomials defined in \eqref{A.68}.
\end{lemma}

\begin{lemma}
Let $k\in\bbZ$. Then the functions $\Phi_+(\cdot,k)|_{\D}$
$($resp., $\Phi_-(\cdot,k)|_{\D}$$)$ are Schur $($resp.,
anti-Schur\,$)$ functions. Moreover, $\Phi_\pm$ satisfy the
following Riccati-type equation
\begin{equation}
\alpha_k
\Phi_\pm(z,k-1)\Phi_\pm(z,k)-\Phi_\pm(z,k-1)+z\Phi_\pm(z,k)
=\ol{\alpha_k}z, \quad z\in\bbC\backslash\dD,\; k\in\Z.
\lb{A.80}
\end{equation}
\end{lemma}

\begin{remark} \lb{rB.17}
$(i)$ In the special case $\alpha=\{\alpha_k\}_{k\in\Z}=0$, one
obtains
\begin{equation}
M_\pm(z,k) = \pm 1, \quad \Phi_+(z,k)=0, \quad
1/\Phi_-(z,k)=0, \quad  z\in\C, \; k\in\Z.
\end{equation}
Thus, strictly speaking, one should always consider
$1/\Phi_-$ rather than $\Phi_-$ and hence refer to the
Riccati-type equation of $1/\Phi_-$,
\begin{equation}
\ol{\alpha_k}z\f{1}{\Phi_-(z,k-1)}\f{1}{\Phi_-(z,k)}
+\f{1}{\Phi_-(z,k)} -z\f{1}{\Phi_-(z,k-1)}=\alpha_k, \quad
z\in\C\backslash\dD, \; k\in\bbZ,  \lb{B.78}
\end{equation}
rather than that of $\Phi_-$, etc. For simplicity of
notation, we will often avoid this distinction between
$\Phi_-$ and $1/\Phi_-$ and usually just invoke $\Phi_-$
whenever confusions are unlikely. \\
$(ii)$ We note that $M_\pm(z,k)$ and $\Phi_\pm(z,k)$,
$z\in\dD$, $k\in\bbZ$, have nontangential limits to $\dD$
$\mu_0$-a.e.\ In particular, the Riccati-type equations
\eqref{A.77}, \eqref{A.80}, and \eqref{B.78} extend to $\dD$
$\mu_0$-a.e.
\end{remark}

\begin{lemma}
Let $z\in\bbC\backslash(\dD\cup\{0\})$ and fix $k_0,
k_1\in\bbZ$. Then  the resolvent $(U-zI)^{-1}$ of the
unitary CMV operator $U$ on $\ell^2(\bbZ)$ is given in
terms of its matrix representation in the standard basis of
$\ltz$ by
\begin{align} \nonumber
&(U-zI)^{-1}(k,k') = \frac{(-1)^{k_1+1}}{z
W\left(\begin{pmatrix} u_+(z,k_1,k_0)\\ v_+(z,k_1,k_0)\end{pmatrix},
\begin{pmatrix}u_-(z,k_1,k_0)\\ v_-(z,k_1,k_0)\end{pmatrix}\right)}  \\
& \quad \times
\begin{cases}
u_-(z,k,k_0)v_+(z,k',k_0), & k < k' \text{ and } k = k'
\text{ odd},
\\
v_-(z,k',k_0) u_+(z,k,k_0), & k' < k \text{ and } k = k'
\text{ even},
\end{cases} \quad k,k' \in\Z, \label{A.84}
\end{align}
where
\begin{align}
& W\left(\binom{u_+(z,k_1,k_0)}{v_+(z,k_1,k_0)},
\binom{u_-(z,k_1,k_0)}{v_-(z,k_1,k_0)}\right) =
\det\left(\begin{pmatrix}
u_+(z,k_1,k_0) & u_-(z,k_1,k_0) \\
v_+(z,k_1,k_0) & v_-(z,k_1,k_0)
\end{pmatrix}\right) \no \\
& \quad = (-1)^{k_1}[M_+(z,k_0)-M_-(z,k_0)]
\begin{cases}
2z, & k_0  \text{ odd}, \\ 2,   & k_0 \text{ even},
\end{cases} \quad k_0,k_1\in\Z. \lb{A.85}
\end{align}
Moreover, since $0\in\bbC\backslash\sigma(U)$, \eqref{A.84}
analytically extends to $z=0$.
\end{lemma}

Next, we briefly turn to Weyl--Titchmarsh theory for CMV
operators on $\bbZ$. We denote by $d\Omega(\cdot,k)$,
$k\in\Z$, the $2 \times 2$ matrix-valued measure,
\begin{align}
d\Omega(\ze,k) &= d
\begin{pmatrix}
\Omega_{0,0}(\ze,k) & \Omega_{0,1}(\ze,k)
\\
\Omega_{1,0}(\ze,k) & \Omega_{1,1}(\ze,k)
\end{pmatrix} \no
\\ &= d
\begin{pmatrix}
(\de_{k-1},E_{U}(\ze)\de_{k-1})_{\ell^2(\Z)}
&(\de_{k-1},E_{U}(\ze)\de_{k})_{\ell^2(\Z)}
\\
(\de_{k},E_{U}(\ze)\de_{k-1})_{\ell^2(\Z)} &
(\de_{k},E_{U}(\ze)\de_{k})_{\ell^2(\Z)}
\end{pmatrix}, \quad \ze \in\dD, \label{A.87}
\end{align}
where $dE_{U}(\cdot)$ denotes the operator-valued spectral
measure of the unitary CVM operator $U$ on $\ell^2(\bbZ)$,
\begin{equation}
U=\oint_{\dD} dE_{U}(\zeta)\,\zeta.
\end{equation}

We also introduce the $2 \times 2$ matrix-valued
function $\cM(\cdot,k)$, $k\in\Z$, by
\begin{align}
&\cM(z,k) = \begin{pmatrix} M_{0,0}(z,k) & M_{0,1}(z,k) \\
M_{1,0}(z,k) & M_{1,1}(z,k) \end{pmatrix} \no
\\ &\quad =
\begin{pmatrix}
(\de_{k-1},(U+zI)(U-zI)^{-1}\de_{k-1})_{\ell^2(\Z)}
&(\de_{k-1},(U+zI)(U-zI)^{-1}\de_{k})_{\ell^2(\Z)}
\\
(\de_{k},(U+zI)(U-zI)^{-1}\de_{k-1})_{\ell^2(\Z)} &
(\de_{k},(U+zI)(U-zI)^{-1}\de_{k})_{\ell^2(\Z)}
\end{pmatrix} \no
\\ &\quad =
\oint_\dD d\Omega(\ze,k)\, \frac{\ze+z}{\ze-z}, \quad
z\in\bbC\backslash\dD. \lb{A.88}
\end{align}
We note that
\begin{align}
M_{0,0}(\cdot,k+1) = M_{1,1}(\cdot,k), \quad k\in\bbZ \lb{A.93}
\end{align}
and
\begin{align}
M_{1,1}(z,k) &= (\de_{k}, (U+zI)(U-zI)^{-1}\de_{k})_{\ell^2(\Z)}
\lb{A.94}
\\ & = \oint_\dD d\Omega_{1,1}(\zeta,k) \, \frac{\zeta+z}{\zeta-z},
\quad z\in\bbC\backslash\dD,\; k\in\Z, \lb{A.95}
\end{align}
where
\begin{equation}
d\Omega_{1,1}(\zeta,k)=d(\de_{k},E_U(\zeta)\de_{k})_{\ell^2(\Z)}, \quad
\zeta\in\dD. \lb{A.96}
\end{equation}
Thus, $M_{0,0}|_{\D}$ and $M_{1,1}|_{\D}$ are Caratheodory
functions. Moreover, by \eqref{A.94} one infers that
\begin{equation}
M_{1,1}(0,k)=1, \quad k\in\Z. \lb{A.98}
\end{equation}

\begin{lemma}
Let $z\in\bbC\backslash\dD$. Then the functions
$M_{\ell,\ell'}(\cdot,k)$, $\ell,\ell'=0,1$, and $M_\pm(\cdot,k)$,
$k\in\bbZ$, satisfy the following relations
\begin{align}
M_{0,0}(z,k) & =1+\f{[\ol{a_k}-\ol{b_k}M_+(z,k)][a_k
+b_kM_-(z,k)]}{\rho_k^2[M_+(z,k)-M_-(z,k)]}, \label{A.97a} \\
M_{1,1}(z,k) &=
\frac{1-M_+(z,k)M_-(z,k)}{M_+(z,k)-M_-(z,k)}, \label{A.97} \\
M_{0,1}(z,k)&=\f{-1}{\rho_k[M_+(z,k)-M_-(z,k)]}
\begin{cases}
{[1-M_+(z,k)][\ol{a_k}-\ol{b_k}M_-(z,k)]}, & k \text{ odd},
\\
{[1+M_+(z,k)][a_k+b_k M_-(z,k)]}, & k \text{ even},
\end{cases} \label{A.97b}
\\
M_{1,0}(z,k)&=\f{-1}{\rho_k[M_+(z,k)-M_-(z,k)]}
\begin{cases}
{[1+M_+(z,k)][a_k+b_kM_-(z,k)]}, & k \text{ odd},
\\
{[1-M_+(z,k)][\ol{a_k}-\ol{b_k}M_-(z,k)]}, & k \text{
even},
\end{cases} \label{A.97c}
\end{align}
where $a_k=1+\al_k$ and $b_k=1-\al_k$, $k\in\Z$.
\end{lemma}

Finally, introducing the functions $\Phi_{1,1}(\cdot,k)$,
$k\in\bbZ$, by
\begin{equation} \label{A.99}
\Phi_{1,1}(z,k) = \frac{M_{1,1}(z,k)-1}{M_{1,1}(z,k)+1}, \quad
z\in\C\backslash\dD,
\end{equation}
then,
\begin{equation}
M_{1,1}(z,k) = \frac{1+\Phi_{1,1}(z,k)}{1-\Phi_{1,1}(z,k)}, \quad
z\in\C\backslash\dD.
\lb{A.100}
\end{equation}

\begin{lemma}
The function $\Phi_{1,1}|_{\D}$ is a Schur function and $\Phi_{1,1}$
is related to $\Phi_\pm$ by
\begin{equation}
\Phi_{1,1}(z,k) = \frac{\Phi_+(z,k)}{\Phi_-(z,k)}, \quad
z\in\C\backslash\dD, \; k\in\bbZ. \lb{A.101}
\end{equation}
\end{lemma}

Denoting by $I_2$ the identity operator in $\bbC^2$, we state the
following result.

\begin{theorem} \lb{tB.22}
Let $k_0\in\bbZ$. Then the CMV operator $U$ on $\ell^2(\Z)$ is
unitarily equivalent to the operator of multiplication by
$I_2 id$ $($where $id(\ze)=\ze$, $\ze\in\dD$$)$
on $L^2(\dD;d\Om(\cdot,k_0))$. Thus,
\begin{equation}
\si(U) = \supp \, (d\Omega(\cdot,k_0))=
\supp \, (d\Omega^{\rm tr}(\cdot,k_0)), \lb{A.102}
\end{equation}
where
\begin{equation}
d\Omega^{\rm tr}(\cdot,k_0) = d\Omega_{0,0}(\cdot,k_0) +
d\Omega_{1,1}(\cdot,k_0) \lb{A.103}
\end{equation}
denotes the trace measure of $d\Omega(\cdot,k_0)$.
\end{theorem}

{\bf Acknowledgments.}
We are indebted to Barry Simon for providing us with a copy of his
forthcoming two-volume treatise \cite{Si04} and for his interest in
this paper.



\begin{thebibliography}{99}
%
\bi{AL75} M.~J.~Ablowitz and J.~F.~Ladik, {\it Nonlinear
differential-difference equations}, J. Math. Phys. {\bf 16}, 598--603
(1975).
%
\bi{APT04} M.~J.~Ablowitz, B.~Prinari, and A.~D.~Trubatch, {\it
{D}iscrete and {C}ontinuous {N}onlinear {S}chr\"odinger {S}ystems},
London Math. Soc. Lecture Note Series, Vol.\ 302, Cambridge
University Press, Cambridge, 2004.
%
\bi{Ak65} N.~I.~Akhiezer, {\it The Classical Moment Problem},
Oliver \& Boyd., Edinburgh, 1965.
%
\bi{AG81} N.~I.~Akhiezer and I.\ M.\ Glazman, {\it Theory of
Operators in Hilbert Space}, Vol. I, Pitman, Boston, 1981.
%
\bi{AD56} N.~Aronszajn and W.~F.~Donoghue, {\it On exponential
representations of analytic functions in the upper half-plane with
positive imaginary part}, J. Anal. Math. {\bf 5}, 321-388 (1956-57).
%
\bibitem{AD64} N.~Aronszajn and W.~F.~Donoghue, {\it A
supplement to the
paper on exponential representations of analytic functions
in the upper half-plane with positive imaginary parts}, J.
Analyse Math. {\bf 12}, 113--127 (1964).
%
\bi{Bo46} G.~Borg, {\it Eine Umkehrung der Sturm-Liouvilleschen
Eigenwertaufgabe}, Acta Math. {\bf 78}, 1--96 (1946).
%
\bi{CMV03} M.\ J.\ Cantero, L.\ Moral, and L.\ Vel\'azquez, {\it
Five-diagonal matrices and zeros of orthogonal polynomials on the
unit circle}, Lin. Algebra Appl. {\bf 362}, 29--56 (2003).
%
\bi{CG02} S.~Clark and F.~Gesztesy, {\it Weyl--Titchmarsh
$M$-function asymptotics and Borg-type theorems for Dirac
operators}, Trans. Amer. Math. Soc. {\bf 354}, 3475--3534 (2002).
%
\bi{CGHL00} S.~Clark, F.~Gesztesy, H.~Holden, and
B.~M.~Levitan, {\it Borg-type theorems for matrix-valued
Schr\"odinger operators}, J. Diff. Eqs. {\bf 167}, 181--210 (2000).
%
\bi{CGR04} S.\ Clark, F.\ Gesztesy, and W.\ Renger, {\it Borg-type
theorems for matrix-valued Jacobi and Dirac finite difference
operators}, preprint, 2004.
%
\bi{Cr89} W.~Craig, {\it The trace formula for Schr\"odinger
operators on the line}, Commun.~Math.~Phys. {\bf 126},
379--407 (1989).
%
\bi{DS83} P.~Deift and B.~Simon, {\it Almost periodic Schr\"odinger
operators III. The absolutely continuous spectrum in one dimension},
Commun.~Math.~Phys. {\bf 90}, 389--411 (1983).
%
\bi{De95} B.~Despr\'es, {\it The Borg theorem for the vectorial Hill's
equation}, Inverse Probl. {\bf 11}, 97--121 (1995).
%
\bi{Fl75} H.\ Flaschka, {\it Discrete and periodic illustrations of
some aspects of the inverse method}, in {\it Dynamical Systems, Theory
and Applications}, J.\ Moser (ed.), Lecture Notes In Physics, Vol.\ 38,
Springer Verlag, Berlin, 1975, p.\ 441--466.
%
\bi{GGH05} J.\ S.\ Geronimo, F.\ Gesztesy, H.\ Holden, {\it
Algebro-geometric solutions of the Baxter--Szeg\H o
difference equation}, to appear in Commun. Math. Phys.
%
\bi{GJ96} J.~S.~Geronimo and R.~Johnson, {\it Rotation number
associated with difference equations satisfied by polynomials
orthogonal on the unit circle}, J. Diff. Eqs. {\bf 132}, 140--178
(1996).
%
\bi{GJ98} J.~S.~Geronimo and R.~Johnson, {\it An inverse problem
associated with polynomials orthogonal on the unit circle}, Commun.
Math. Phys. {\bf 193}, 125--150 (1998).
%
\bi{GT94} J.~S.~Geronimo and A.~Teplyaev, {\it A difference equation
arising from the trigonometric moment problem having random reflection
coefficients--an operator theoretic approach}, J. Funct. Anal.
{\bf 123}, 12--45 (1994).
%
\bi{Ge46} J.~Geronimus, {\it On the trigonometric moment problem},
Ann. Math. {\bf 47}, 742--761 (1946).
%
\bi{Ge48} Ya.~L.~Geronimus, {\it Polynomials orthogonal on a circle and
their applications}, Commun. Soc. Mat. Kharkov {\bf 15}, 35--120 (1948);
Amer. Math. Soc. Transl. (1) {\bf 3}, 1--78 (1962).
%
\bi{Ge61} Ya.~L.~Geronimus, {\it {O}rthogonal {P}olynomials},
Consultants Bureau, New York, 1961.
%
\bi{GH05} F.~Gesztesy and H.~Holden, {\it {S}oliton {E}quations and
{T}heir {A}lgebro-{G}eometric {S}olutions. Volume {II}:
$(1+1)$-Dimensional Discrete Models}, Cambridge Studies in Adv.
Math., Cambridge University Press, Cambridge, in preparation.
%
\bi{GKT96} F.~Gesztesy, M.~Krishna, and G.~Teschl, {\it On
isospectral sets of Jacobi operators,} Commun. Math. Phys.
{\bf 181}, 631--645 (1996).
%
\bi{GS96} F.~Gesztesy and B.~Simon, {\it The $\xi$ function},
Acta Math. {\bf 176}, 49--71 (1996).
%
\bi{GT00} F.~Gesztesy and E.~Tsekanovskii, {\it On
matrix-valued Herglotz functions}, Math. Nachr. {\bf 218}, 61--138
(2000).
%
\bi{GZ05} F.\ Gesztesy and M.\ Zinchenko, {\it Weyl--Titchmarsh
theory for CMV operators associated with orthogonal polynmials on
the unit circle}, preprint, 2004.
%
\bi{GJ84} R.~Giachetti and R.~A.~Johnson, {\it Spectral
theory of second-order almost periodic differential operators
and its relation to classes of nonlinear evolution equations},
Nuovo Cim. {\bf 82B}, 125--168 (1984).
%
\bi{GJ86} R.~Giachetti and R.~A.~Johnson, {\it The Floquet exponent
for two-dimensional linear systems with bounded coefficients}, J.
Math. pures et appl. {\bf 65}, 93--117 (1986).
%
\bi{GN01} L.~Golinskii and P.~Nevai, {\it {S}zeg\H o difference
equations, transfer matrices and orthogonal polynomials on the unit
circle}, Commun. Math. Phys. {\bf 223}, 223--259 (2001).
%
\bi{Gr60} D.~S.~Greenstein, {\it On the analytic continuation of
functions which map the upper half plane into itself}, J. Math.
Anal. Appl. {\bf 1}, 355--362 (1960).
%
\bi{GS84}
U.~Grenander and G.~Szeg\H o, {\it {T}oeplitz {F}orms and their
{A}pplications}, University of California Press, Berkeley, 1958; 2nd
ed., Chelsea, New York, 1984.
%
\bi{Jo82} R.~A.~Johnson, {\it The recurrent {H}ill's equation},
J. Diff. Eqs. {\bf 46}, 165--193 (1982).
%
\bibitem{KK74} I.~S.~Kac and M.~G.~Krein, {\it $R$-functions--analytic
functions mapping the upper halfplane into itself}, Amer. Math. Soc.
Transl. (2) {\bf 103}, 1-18 (1974).
%
\bi{Ko84} S.~Kotani, {\it Ljapunov indices determine absolutely
continuous spectra of stationary random one-dimensional
Schr\"odinger operators}, in {\it Stochastic Analysis\/}, K.~It{\v o}
(ed.), North-Holland, Amsterdam, 1984, pp.\ 225--247.
%
\bi{Ko87} S.~Kotani, {\it One-dimensional random Schr\"odinger operators
and Herglotz functions}, in {\it Probabilistic Methods
in Mathematical Physics\/}, K.~It{\v o} and N.~Ikeda (eds.), Academic
Press, New York, 1987, pp.~219--250.
%
\bi{KK88} S.~Kotani and M.~Krishna, {\it Almost periodicity of some
random potentials}, J.~Funct.~Anal. {\bf 78}, 390--405 (1988).
%
\bi{Kr45}
M.~G.~Krein, {\it On a generalization of some investigations
of {G}.~{S}zeg\H o, {V}.~{S}mirnoff, and {A}.~{K}olmogoroff},
Doklady Akad. Nauk SSSR {\bf 46}, 91--94 (1945). ({R}ussian).
%
\bi{Ma94} M.~M.~Malamud, {\it Similarity of Volterra
operators and
related questions of the theory of differential equations
of fractional
order}, Trans.~Moscow Math.~Soc. {\bf 55}, 57--122 (1994).
%
\bi{Ma99a} M.~M.~Malamud,
{\it Borg type theorems for first-order systems on a finite
interval}, Funct. Anal. Appl. {\bf 33}, 64--68 (1999).
%
\bi{MEKL95} P.~D.~Miller, N.~M.~Ercolani, I.~M.~Krichever, and
C.~D.~Levermore, {\it {F}inite genus solutions to the
{A}blowitz--{L}adik equations}, Comm. Pure Appl.  Math.
{\bf 4}, 1369--1440 (1995).
%
\bi{NS05} I.~Nenciu, {\it Lax pairs for the Ablowitz-Ladik system via
orthogonal polynomials on the unit circle}, Int. Math. Res. Notices, to 
appear.
%
\bi{PY04} F.~Peherstorfer and P.~Yuditskii, {\it Asymptotic behavior
of polynomials orthonormal on a homogeneous set}, J. Analyse Math.
{\bf 89}, 113--154 (2003).
%
\bi{Si04a} B.~Simon, {\it Analogs of the $m$-function in the theory of
orthogonal polynomials on the unit circle}, J. Comp. Appl. Math.
{\bf 171}, 411--424 (2004).
%
\bi{Si04}
B.~Simon, {\it Orthogonal Polynomials on the Unit Circle, Part 1: 
Classical Theory, Part 2: Spectral Theory}, AMS Colloquium Publication 
Series, Vol.\ 54, Providence, R.I., 2005.
%
\bi{Si04b} B.~Simon, {\it Orthogonal polynomials on the unit circle:
New results}, {\it Intl.\ Math.\ Res.\ Notices}, 2004, No.\ 53,
2837--2880.
%
\bi{SY95} M.~Sodin and P.~Yuditskii, {\it Almost periodic Sturm-Liouville
operators with Cantor  homogeneous spectrum}, Comment. Math. Helvetici
{\bf 70}, 639--658 (1995).
%
\bi{SY96a} M.~Sodin and P.~Yuditskii, Almost periodic Jacobi matrices
with homogeneous spectrum, infinite dimensional Jacobi inversion, and
Hardy spaces of character-automorphic functions, J. Geom. Anal. {\bf
7}, 387--435 (1997).
%
\bi{Sz20} G.~Szeg{\H o}, {\it {B}eitr\"{a}ge zur {T}heorie der
{T}oeplitzschen {F}ormen I}, Math. Z. {\bf 6}, 167--202 (1920).
%
\bi{Sz21} G.~Szeg{\H o}, {\it {B}eitr\"{a}ge zur {T}heorie der
{T}oeplitzschen {F}ormen II}, Math. Z. {\bf 9}, 167--190 (1921).
%
\bi{Sz78} G.~Szeg{\H o}, {\it {O}rthogonal {P}olynomials}, Amer Math.
Soc. Colloq. Publ., Vol.\ 23, Amer. Math. Soc., Providence, R.I.,
1978.
%
\bi{Te98} G.~Teschl, {\it Trace formulas and inverse
spectral theory for Jacobi operators}, Commun. Math. Phys.
{\bf 196}, 175--202 (1998).
%
\bi{Te00} G.~Teschl, {\it Jacobi Operators and Completely Integrable
Nonlinear Lattices}, Mathematical Surveys and Monographs, {\bf 72},
American Mathematical Society, (2000).
%
\bi{To63} Ju.\ Ja.\ Tom{\v c}uk, {\it Orthogonal polynomials on a given
system of arcs of the unit circle}, Sov. Math. Dokl. {\bf 4}, 931--934
(1963).
%
\bi{Ve35} S.\ Verblunsky, {\it On positive harmonic functions: A
contribution to the algebra  of {F}ourier series},
Proc. London Math. Soc. (2) {\bf 38}, 125--157 (1935).
%
\bi{Ve36} S.\ Verblunsky, {\it On positive harmonic functions (second
paper)}, Proc. London Math. Soc. (2) {\bf 40}, 290--320 (1936).
%
\end{thebibliography}
\end{document}